# TIGHTNESS FOR A FAMILY OF RECURSION EQUATIONS

By Maury Bramson[1] and Ofer Zeitouni[2]

*University of Minnesota, and University of Minnesota and Weizmann Institute of Science*

In this paper we study the tightness of solutions for a family of recursion equations. These equations arise naturally in the study of random walks on tree-like structures. Examples include the maximal displacement of a branching random walk in one dimension and the cover time of a symmetric simple random walk on regular binary trees. Recursion equations associated with the distribution functions of these quantities have been used to establish weak laws of large numbers. Here, we use these recursion equations to establish the tightness of the corresponding sequences of distribution functions after appropriate centering. We phrase our results in a fairly general context, which we hope will facilitate their application in other settings.

**1. Introduction.** Branching random walks (BRW) have been studied since the 1960's, with various cases having been studied earlier (see, e.g., [6] and [19]). As the name suggests, such processes consist of particles that execute random walks while also branching. One typically assumes that the particles lie on $\mathbb{R}$, and that the corresponding discrete time process starts from a single particle at 0. When the branching is supercritical, that is, particles on the average have more than one offspring, the number of particles will grow geometrically in time off the set of extinction.

Two types of results have received considerable attention. The first pertains to the location of the main body of particles, and states roughly that, at large times, this distribution is given by a suitably scaled normal distribution, if typical tail conditions on the random walk increments and on the offspring distribution are assumed (see, e.g., [8] and [20]). The second pertains to the maximal displacement of BRW, that is, the position of the

Received February 2007; revised March 2008.
[1]Supported in part by NSF Grant DMS-02-26245 and CCF-0729537.
[2]Supported in part by NSF Grants DMS-03-02230 and DMS-05-03775.
*AMS 2000 subject classifications.* 60J80, 60G50, 39B12.
*Key words and phrases.* Tightness, recursion equations, branching random walk, cover time.

particle furthest to the right. Results date back to [18], who demonstrated a weak law of large numbers. There have since been certain refinements for the limiting behavior of the maximal displacement, but the theory remains incomplete. Here, we obtain sufficient conditions for the tightness of the maximal displacement after centering. Our result is motivated by analogous behavior for the maximal displacement of branching Brownian motion (BBM), about which much more is known.

Another problem concerns the cover time for regular binary trees $\mathbf{T}_n$ of depth $n$. A particle, which starts at the root of the tree, executes a symmetric simple random walk on $\mathbf{T}_n$. How long does it take for the particle to visit every site in $\mathbf{T}_n$? In [4], a weak law of large numbers was given for regular $k$-ary trees as $n \to \infty$. (That is, each parent has precisely $k$ offspring.) Here, we show that, under an appropriate scaling, the sequence of square roots of the cover times of $\mathbf{T}_n$ is tight and nondegenerate, after centering. (The same result also holds for regular $k$-ary trees.)

Distribution functions associated with the maximal displacement of BRW and cover times for trees are known to satisfy recursion equations. (See [5] for a survey of recursions arising in these and similar contexts.) The distribution function itself is used in the context of the maximal displacement; for cover times, the relationship is more complicated. These recursion equations were used in [18] and [4] for the weak laws of large numbers they established. We will employ these equations to demonstrate our results on tightness. We will phrase our main result on tightness in a more general setting, which includes both the maximal displacement of BRW and the cover time of trees.

We next describe the maximal displacement and cover time problems in detail and state the corresponding results, Theorems 1.1 and 1.2. We then summarize the remainder of the paper. We will postpone until Section 2 our general result on tightness, Theorem 2.5, because of the terminology that is required.

We define branching random walk on $\mathbb{R}$ formally as follows. A particle starting at 0 moves randomly to a site according to a given distribution function $G(\cdot)$, at which time it dies and gives birth to $k$ offspring with probability $p_k$, independently of the previous motion. Each of these offspring, in turn, moves independently according to the same distribution $G(\cdot)$ over the next time step, then dies and gives birth to $k$ offspring according to the distribution $\{p_k\}$. This procedure is repeated at integer times, with the movement of all particles and the number of offspring being assumed to be independent of one another. To avoid the possibility of extinction and trivial special cases, we assume that $p_0 = 0$ and $p_1 < 1$. This implies that the mean number of offspring $m_1 = \sum_{k=1}^{\infty} k p_k > 1$.

Let $Z_n$ denote the number of particles at time $n$ of the BRW, with $\mathfrak{x}_k(n), k = 1, \ldots, Z_n$, being the positions of these particles when ordered in

some fashion. We write

$$\mathcal{M}_n = \max_{1 \leq k \leq Z_n} \mathfrak{x}_k(n) \tag{1.1}$$

for the *maximal displacement* of the BRW at time $n$. [The *minimal displacement* will be given by $\mathcal{M}_n^{\min} = \min_{1 \leq k \leq Z_n} \mathfrak{x}_k(n)$. Questions about $\mathcal{M}_n$ or $\mathcal{M}_n^{\min}$ can be rephrased as questions about the other by substituting $-x$ for the coordinate $x$ and reflecting $G(\cdot)$ accordingly.] When $G(0) = 0$, one can alternatively interpret $G(\cdot)$ as the lifetime distribution of individual particles of the branching process. In this setting, $\mathcal{M}_{n+1}$ becomes the last death time of the $n$th generation of particles. ($\mathcal{M}_{n+1}^{\min}$ then becomes the first death time.)

Let $F_n(\cdot)$ denote the distribution function of $\mathcal{M}_n$, and set $\bar{F}_n(\cdot) = 1 - F_n(\cdot)$. One can express $\bar{F}_n(\cdot)$ recursively in terms of $\bar{G}(\cdot) = 1 - G(\cdot)$ and

$$Q(u) = 1 - \sum_k p_k (1-u)^k \qquad \text{for } u \in [0,1]. \tag{1.2}$$

One has

$$\bar{F}_{n+1}(x) = -(\bar{G} * Q(\bar{F}_n))(x) = -\int_{y \in \mathbb{R}} \bar{G}(x-y) \, dQ(\bar{F}_n(y)), \tag{1.3}$$

with $\bar{F}_0(x) = \mathbf{1}_{\{x<0\}}$. Here, $*$ denotes the standard convolution. [One requires the minus sign since the function $\bar{F}_n(\cdot)$ is decreasing.] Equation (1.3) is the backward equation for $\bar{F}_{n+1}(\cdot)$ in terms of $\bar{F}_n(\cdot)$. It is simple to derive by relating it to the maximal displacement of the $n$th generation offspring for each of the first generation offspring of the original particle. The composite function $Q(\bar{F}_n)$ gives the distribution of the maximum of these individual maximal displacements (relative to their parents), with convolution by $\bar{G}$ then contributing the common displacement due to the movement of the original particle. In the special case where the branching is binary, that is, where $p_2 = 1$, (1.2) reduces to $Q(u) = 2u - u^2$. We note that $Q: [0,1] \to [0,1]$ is strictly concave, in general, with

$$Q(0) = 0, \qquad Q(1) = 1 \quad \text{and} \quad Q'(0) = m_1 > 1. \tag{1.4}$$

One can equally well assume that branching for the BRW occurs at the beginning of each time step, before the particles move rather than after. The corresponding distribution functions $F_n^r(\cdot)$ then satisfy the analog of (1.3),

$$\bar{F}_{n+1}^r = Q(-\bar{G} * \bar{F}_n^r). \tag{1.5}$$

Since $\bar{F}_1 = -\bar{G} * \bar{F}_0$, one has $\bar{F}_n^r = Q(\bar{F}_n)$ for all $n$; consequently, $\{F_n\}$ and $\{F_n^r\}$ will have the same asymptotic behavior. The distribution functions $F_n^{r,\min}(\cdot)$ of the minimal displacement of this BRW were studied in [18]. They satisfy a recursion equation that is similar by (1.5).

It follows from [18], Theorem 2, that for appropriate $\gamma_0$,

(1.6) $$F_n^{r,\min}(\gamma n) \to \begin{cases} 0, & \text{for } \gamma < \gamma_0, \\ 1, & \text{for } \gamma > \gamma_0, \end{cases}$$

as $n \to \infty$, provided $G(\cdot)$ has finite mean and its support is bounded below. [Related results were proved in [21] and [22], and by H. Kesten (unpublished).] Hammersley believed that the minimal displacement $\mathcal{M}_n^{r,\min}$ was in some sense subadditive. This was made precise in [24]; the subadditive ergodic theorem given there demonstrates the strong law analog of (1.6). (The strong law was demonstrated using other techniques in [21].) Analogous laws of large numbers hold for $F_n(\cdot)$, $F_n^{\min}(\cdot)$ and $F_n^r(\cdot)$, that is, for $\mathcal{M}_n, \mathcal{M}_n^{\min}$ and $\mathcal{M}_n^r$. In this paper we will investigate the refined behavior of $F_n(\cdot)$.

There is an older, related theory of branching Brownian motion. Individual particles are assumed to execute standard Brownian motion on $\mathbb{R}$. Starting with a single particle at 0, particles die after independent rate-1 exponentially distributed holding times, at which point they give birth to $k$ offspring with distribution $\{p_k\}_{k \geq 1}$. All particles are assumed to move independently of one another and of the number of offspring at different times, which are themselves independent. The maximal displacement

$$\mathcal{M}_t = \max_{1 \leq k \leq Z_t} \mathfrak{r}_k(t)$$

is the analog of (1.1), where, as before, $Z_t$ and $\mathfrak{r}_k(t)$, $k = 1, \ldots, Z_t$, are the number of particles and their positions at time $t$. It is not difficult to show (see, e.g., [28]) that $u(t, x) = P(\mathcal{M}_t > x)$ satisfies

(1.7) $$u_t = \tfrac{1}{2} u_{xx} + f(u),$$

with

(1.8) $$f(u) = Q(u) - u$$

and $u(0, x) = \mathbf{1}_{\{x < 0\}}$. When the branching is binary, $f(u) = u(1 - u)$.

When $f(\cdot)$ is continuously differentiable and satisfies the more general equation

(1.9) $\quad f(0) = f(1) = 0, \qquad f(u) > 0, \qquad f'(u) \leq f'(0), \qquad \text{for } u \in (0, 1),$

(1.7) is typically either referred to as the *K–P–P equation* or the *Fisher equation*. For solutions $u(t, x)$ of (1.7) with $u(0, x) = \mathbf{1}_{\{x < 0\}}$, $1 - u(t, \cdot)$ will be a distribution function for each $t$. In both [23] and [17], (1.7) was employed to model the spread of an advantageous gene through a population.

In [23], it was shown that, under (1.9) and $u(0, x) = \mathbf{1}_{\{x < 0\}}$, the solution of (1.7) converges to a traveling wave $w(x)$, in the sense that, for appropriate $m(t)$,

(1.10) $$u(t, x + m(t)) \to w(x) \qquad \text{as } t \to \infty$$

uniformly in $x$, where $1 - w(\cdot)$ is a distribution function for which $\tilde{u}(t,x) = w(x - \sqrt{2}t)$ satisfies (1.7). Moreover,

$$(1.11) \qquad m(t)/t \to \sqrt{2} \qquad \text{as } t \to \infty.$$

[The centering term $m(t)$ can be chosen so that $u(t, m(t)) = c$, for any given $c \in (0,1)$.] In particular,

$$u(t, \gamma t) \to \begin{cases} 1, & \text{for } \gamma < \sqrt{2}, \\ 0, & \text{for } \gamma > \sqrt{2}, \end{cases}$$

which is the analog of (1.6).

A crucial step in the proof of (1.10) consists of showing that

$$(1.12) \quad \begin{aligned} &u(t, x + m(t)) \text{ is decreasing in } t \text{ for } x < 0, \\ &u(t, x + m(t)) \text{ is increasing in } t \text{ for } x > 0. \end{aligned}$$

That is, $v(t,\cdot) = u(t, \cdot + m(t))$ "stretches" as $t$ increases. A direct consequence of (1.10) and (1.12) is that the family $v(t,\cdot)$ is *tight*, that is, for each $\varepsilon > 0$, there is an $A_\varepsilon > 0$, so that, for all $t$,

$$(1.13) \qquad v(t, -A_\varepsilon) - v(t, A_\varepsilon) > 1 - \varepsilon.$$

More can be said about $m(\cdot)$ and the convergence of $v(t,\cdot)$ under more general initial data ([9] and [10]).

Although BRW is the discrete time analog of branching Brownian motion, with (1.3) corresponding to (1.7), more refined results on the asymptotic behavior of $\bar{F}_n(\cdot)$ corresponding to those of $u(t,\cdot)$ in (1.10) have, except in special cases, remained elusive. When $G(\cdot)$ admits a density which is logarithmically concave, that is, $G(\cdot)$ satisfies

$$(1.14) \qquad G'(x) = e^{-\varphi(x)}, \qquad \text{where } \varphi(x) \in (-\infty, \infty] \text{ is convex},$$

one can show that the analog of (1.12) holds for $\bar{F}_0(x) = \mathbf{1}_{\{x<0\}}$. As in [26] and [7], the analog of (1.10) follows from this. Results of this nature for general $G(\cdot)$ are not known. In fact, without some modification, the analog of (1.10) will be false in general, as when $G(\cdot)$ is concentrated on the integers and $\gamma_0 \notin \mathbb{Z}$.

There has recently also been some interest in related problems that arise in the context of sorting algorithms, for which the movement of offspring of a common parent will be dependent. (BRW with such dependence are also well known in the general branching literature; see, e.g., [13, 21, 27].) [15] showed the analog of (1.10) for a specific at choice of $G(\cdot)$. In [29] and in [2] [in the latter paper, for general $G(\cdot)$ having bounded support], $m(t)$ is calculated for related models. [12] treats a generalization of the model in [15].

In this paper we will show that, after appropriate centering, the sequence $\{\mathcal{M}_n\}_{n\geq 0}$ corresponding to the maximal displacement of BRW is tight. The shifted sequence $\{\mathcal{M}_n^s\}_{n\geq 0}$ is given by

$$\mathcal{M}_n^s = \mathcal{M}_n - \text{Med}(F_n), \tag{1.15}$$

where $\text{Med}(F_n) = \inf\{x : F_n(x) \geq 1/2\}$ and $F_n(\cdot)$ is the distribution function of $\mathcal{M}_n$. $F_n^s(\cdot)$ denotes the distribution function of $\mathcal{M}_n^s$. The sequence $\{\mathcal{M}_n^s\}_{n\geq 0}$ or, equivalently, $\{F_n^s\}_{n\geq 0}$, is *tight* if for any $\varepsilon > 0$, there is an $A_\varepsilon > 0$ such that $F_n^s(A_\varepsilon) - F_n^s(-A_\varepsilon) > 1 - \varepsilon$ for all $n$; this is the analog of (1.13).

Rather than (1.14) as our main condition, we assume that, for some $a > 0$ and $M_0 > 0$, $\bar{G}(\cdot)$ satisfies

$$\bar{G}(x + M) \leq e^{-aM}\bar{G}(x) \qquad \text{for all } x \geq 0, M \geq M_0. \tag{1.16}$$

In addition to specifying that $\bar{G}(\cdot)$ has an exponentially decreasing right tail, (1.16) requires that $\bar{G}(\cdot)$ be "flat" on no interval $[x, x + M]$, for $x$ and $M$ chosen as above. It follows with a little work from [18] (3.97), that in order for $\gamma_0 < \infty$ to hold, the right tail of $\bar{G}(\cdot)$ needs to be exponentially decreasing. The flatness condition included in (1.16) is needed for our method of proof; this additional condition will be satisfied for most distributions that one encounters in practice. We will also require that the branching law for the BRW satisfy $p_1 < 1$ and

$$\sum_{k=1}^{\infty} k^\theta p_k = m_\theta < \infty \qquad \text{for some } \theta \in (1, 2]. \tag{1.17}$$

Employing the above conditions, we now state our main result for branching random walks:

THEOREM 1.1. *Assume that the random walk increments $G(\cdot)$ of a BRW satisfy (1.16) and that the branching law $\{p_k\}_{k\geq 1}$ satisfies $p_1 < 1$ and (1.17). Then, the sequence of random variables $\{\mathcal{M}_n^s\}_{n\geq 0}$ is tight.*

As indicated earlier in the section, Theorem 1.1 will follow from a more general result, Theorem 2.5, whose statement is postponed until the next section. Theorem 2.5 assumes a more general version of the recursion equation (1.3). The heart of the paper, Sections 2 and 3, is devoted to demonstrating Theorem 2.5. Since the distribution function of the maximal displacement of BRW satisfies (1.3), Theorem 1.1 will follow quickly. This is shown in the first part of Section 4.

Our other problem concerns the cover time for regular binary trees. The regular binary tree $\mathbf{T}_n$ of depth $n$ consists of the first $n$ generations, or *levels*, of a regular binary tree, with the root $o$ denoting the original ancestor and $\mathbf{L}_k$ consisting of the set of $2^k$ vertices at the $k$th level, for $k \leq n$, of the corresponding descendents. We consider each level $k-1$ vertex to be an immediate neighbor of the two level $k$ vertices that are immediately descended from it.

In this setting $\{X_j\}_{j \geq 0}$ denotes a symmetric nearest neighbor random walk on $\mathbf{T}_n$, with $X_0 = o$, and with each neighbor being chosen with equal probability at each time step. The *cover time* of $\mathbf{T}_n$ is given by

$$\mathcal{C}_n = \min\left\{J \geq 0 : \bigcup_{j=0}^{J} \{X_j\} = \mathbf{T}_n\right\}.$$

In [4] it was shown that

(1.18) $\qquad \mathcal{C}_n/4(\log 2)n^2 2^n \to_{n \to \infty} 1 \qquad$ in probability.

A natural question is how $\mathcal{C}_n$ should be scaled so that the resulting random variables, after shifting by their medians, are tight. It turns out that the correct scaling is given by

(1.19) $$\mathcal{E}_n = \sqrt{\mathcal{C}_n/2^n}.$$

Defining the shift $\mathcal{E}_n^s = \mathcal{E}_n - \text{Med}(\mathcal{E}_n)$ similarly to (1.15), we will show the following result:

THEOREM 1.2. *The sequence of random variables $\{\mathcal{E}_n^s\}_{n \geq 0}$ for the regular binary tree is tight. Furthermore, it is nondegenerate in the sense that there exists a constant $V > 0$ such that*

(1.20) $$\limsup_{n \to \infty} P(|\mathcal{E}_n^s| < V) < 1.$$

Theorem 1.2 will also follow as a special case of Theorem 2.5; this is shown in the second part of Section 4.2. Considerably more work is required in this case than was required for showing Theorem 1.1. In particular, an intermediate tightness result first needs to be demonstrated in Theorem 4.3, for a random sequence $\{\beta_n\}_{n \geq 0}$ that is introduced in (4.5).

In the short Section 5 we will mention a conjecture on the tightness of the cover times for the lattice tori $\mathbb{Z}_n^2 = \mathbb{Z}^2/n\mathbb{Z}^2$. The intuition behind its proof should rely heavily on the proof of Theorem 1.2.

Theorem 2.5 is phrased quite generally so as to allow it to be applied in other settings. In our way of thinking, the tightness exhibited in Theorems 1.1 and 1.2 are examples of a fairly general phenomenon. Because of this level of generality, the assumptions in Theorem 2.5 require some preparation, which is done at the beginning of Section 2. For a first reading, one may think of the quantities introduced there as generalizations of the branching law $\{p_k\}_{k\geq 1}$ and of the distribution of the random walk increments $G(\cdot)$ employed for BRW.

The proof of Theorem 2.5 breaks into two main steps. The tightness of the distribution functions $\{F_n^s\}_{n\geq 0}$ given there can be shown by showing the corresponding tightness of (a) the right tail of $\{\bar{F}_n^s\}_{n\geq 0}$ and (b) the left tail of $\{\bar{F}_n^s\}_{n\geq 0}$. An important ingredient is the Lyapunov function $L(\cdot)$ that is introduced in (2.12). Section 3 is devoted to showing the bound (2.14) on $\sup_n L(\bar{F}_n)$, in Theorem 2.7, that is needed to show the tightness of the right tail.

The argument for (b) is comparatively quick and is done in the latter part of Section 2. The idea, in spirit, is to show that $\bar{F}_n(\cdot)$ must grow rapidly through successive iterations until reaching values close to 1, at coordinates not changing by much after an individual iteration. One then employs the resulting bound, together with the tightness of the right tail from (a), to obtain the tightness of the left tail. This is done in Proposition 2.9 and Lemma 2.10.

A summary of the results of this paper, without proofs, is given in [11]. There, we provide intuition for the basic steps of our reasoning in the setting of branching random walks.

**2. Definitions and statement of main result.** We begin the section with definitions that are employed in the statement of Theorem 2.5, which is our main result on tightness. After stating Theorem 2.5, we introduce the Lyapunov function that will be our main tool in proving the theorem. As indicated in the introduction, we will then demonstrate the easier half of the theorem, with the more difficult half being demonstrated in Section 3.

We denote by $\mathcal{D}$ the set of functions $\bar{F}: \mathbb{R} \to [0,1]$, such that $F = 1 - \bar{F}$ is a distribution function. We will study sequences of distribution functions $\{F_n\}_{n\geq 0}$ that solve recursions of the form

$$\bar{F}_{n+1} = T_n \bar{F}_n, \tag{2.1}$$

with $\bar{F}_n = 1 - F_n$, where $T_n: \mathcal{D} \to \mathcal{D}$ is chosen from a family of operators that will be defined shortly. Equations of the form (2.1) include the recursions mentioned in the introduction in the context of Theorems 1.1 and 1.2. As mentioned there, we are interested in proving the tightness of $F_n$ after recentering $F_n$ around its median.

We next describe the class of operators $T_n$ for which we carry through our program. For $u, \overline{G}_n(y, \cdot) \in \mathcal{D}$, with $n \in \mathbb{N}$ and $y \in \mathbb{R}$, set

$$(\overline{G}_n \circledast u)(x) = -\int_{y \in \mathbb{R}} \overline{G}_n(y, x-y) \, du(y). \tag{2.2}$$

One may think of $\{\overline{G}_n\}$ as corresponding to a family of random variables $\{N_{y,n}\}_{y\in\mathbb{R}, n\in\mathbb{N}}$, with $P(N_{y,n} > x) = \overline{G}_n(y, x)$ and $u$ as corresponding to an independent random variable $Y$, with $P(Y > y) = u(y)$; (2.2) is then equivalent to

$$(\overline{G}_n \circledast u)(x) = P(Y + N_{Y,n} > x).$$

When $\overline{G}_n(y,x) = \overline{G}_n(x)$, the $\circledast$ operation reduces to the standard convolution, up to a minus sign.

Let $\{Q_n\}$ be a sequence of increasing, nonlinear functions mapping $[0,1]$ onto itself. Our main result concerns recursions of the form

$$(2.3) \quad (T_n u)(x) = (\overline{G}_n \circledast (Q_n(u)))(x) = -\int_{y \in \mathbb{R}} \overline{G}_n(y, x-y) \, dQ_n(u(y)),$$

for $u \in \mathcal{D}$. We first state our assumptions on $\{Q_n\}$ and $\{\overline{G}_n\}$.

Assumption 2.1 gives conditions on the growth of $Q_n$. The first part of the assumption includes a lower bound for the growth of $Q_n$ near 0; the second part gives an upper bound on the concavity of $Q_n$ near 0. A prototypical example is given by $Q_n(x) = 2x - x^2$ for all $n$. (Throughout the following definitions, all constants will be independent of $n$ unless explicitly mentioned otherwise.)

ASSUMPTION 2.1. The functions $Q_n : [0,1] \to [0,1]$ are increasing with $Q_n(0) = 0$ and $Q_n(1) = 1$. Moreover, there exist constants $\delta_0 > 0$, $m \in (1,2]$, $c^* > 0$ and $\theta^* > 0$, such that

$$(2.4) \qquad Q_n(x) \geq mx \qquad \text{for all } x \leq 2\delta_0$$

and

$$(2.5) \quad \frac{x_2}{x_1} \leq \frac{Q_n(x_2)}{Q_n(x_1)}[1 + c^*(Q_n(x_1))^{\theta^*}] \qquad \text{for all } 0 < x_1 < x_2 \leq 2(\delta_0 \wedge x_1).$$

The next assumption gives conditions on the "convolution" kernels $\overline{G}_n$. The condition (G1) specifies monotonicity requirements for $\overline{G}_n$, (G2) specifies a particular exponential tail condition, and (G3) is a centering condition.

ASSUMPTION 2.2. (G1) The functions $\overline{G}_n(y, x-y)$ are increasing in $y$, whereas $\overline{G}_n(y, x)$ are decreasing in $y$.

(G2) There exist constants $a \in (0,1)$ and $M_0 > 0$ such that, for all $x \geq 0$, $y \in \mathbb{R}$ and $M \geq M_0$,

$$(2.6) \qquad \overline{G}_n(y - M, x + M) \leq e^{-aM} \overline{G}_n(y, x).$$

(G3) Choosing $m$ as in Assumption 2.1 and setting $\varepsilon_0 = (\log m)/100$,

$$(2.7) \qquad \overline{G}_n(y, 0) \geq 1 - \varepsilon_0 \qquad \text{for all } y.$$

In applications, the particular choice of $\varepsilon_0$ in (2.7) can often be obtained from $\varepsilon_0 > 0$ by applying an appropriate translation. In this context, it is useful to note that (G2) follows from (G1) and the following superficially weaker assumption:

(G2$'$) There exist constants $a' \in (0,1)$, $L > 0$ and $M' \geq 2L$ such that, for all $x \geq L$ and $y \in \mathbb{R}$,

$$\overline{G}_n(y - M', x + M') \leq e^{-a'M'}\overline{G}_n(y, x). \tag{2.8}$$

One can check that (G2) follows from (G1) and (G2$'$) by setting $a = a'/3$ and $M_0 = 2M'$; this follows from iterating the inequality

$$\overline{G}_n(y - M, x + M) \leq e^{-a'M'}\overline{G}_n(y - M + M', x + M - M')$$

$k - 1$ times, and then applying (G1).

Assumption 2.4 below gives uniform regularity conditions on the transformations $T_n$ defined in (2.3). It gives upper and lower bounds on $T_n u$ in terms of an appropriate nonlinear function $\tilde{Q}$ of translates of $u$. Such a $\tilde{Q}$ will be required to belong to the set of functions $\widetilde{\mathcal{Q}}$ satisfying the following properties.

DEFINITION 2.3. $\widetilde{\mathcal{Q}}$ is the collection of strictly increasing continuous functions $\tilde{Q}\colon [0,1] \mapsto [0,1]$ with $\tilde{Q}(0) = 0, \tilde{Q}(1) = 1$, such that:

(T1) $\tilde{Q}(x) > x$ for all $x \in (0,1)$, and for any $\delta > 0$, there exists $c_\delta > 1$ such that $\tilde{Q}(x) > c_\delta x$ for $x \leq 1 - \delta$.

(T2) For each $\delta \in (0,1)$, there exists a nonnegative function $g_\delta(\varepsilon) \to 0$ as $\varepsilon \searrow 0$, such that if $x \geq \delta$ and $\tilde{Q}(x) \leq (1-\delta)/(1+\varepsilon)$, then

$$\tilde{Q}((1 + g_\delta(\varepsilon))x) \geq (1 + \varepsilon)\tilde{Q}(x). \tag{2.9}$$

(T1) of the preceding definition gives lower bounds on the linear growth of $\tilde{Q} \in \widetilde{\mathcal{Q}}$ away from 1, and (T2) gives a uniform local lower bound on the growth of $\tilde{Q}$ away from 0 and 1. In our applications in Section 4 we will have $Q_n = Q$, with $Q$ strictly concave, in which case we will set $\tilde{Q} = Q$. (T1) and (T2) will then be automatically satisfied. Note that the conditions (2.4) and (2.5) in Assumption 2.1 specify the behavior of $Q_n$ near 0, whereas (T1) and (T2) specify the behavior of $\tilde{Q}$ over all of $(0,1)$.

Employing an appropriate $\tilde{Q} \in \widetilde{\mathcal{Q}}$, we now state our last assumption. It will only be needed for the proof of Lemma 2.10.

ASSUMPTION 2.4. There exists a $\tilde{Q} \in \widetilde{\mathcal{Q}}$ satisfying the following. For each $\eta_1 > 0$, there exists $B = B(\eta_1) \geq 0$ satisfying, for all $u \in \mathcal{D}$ and $x \in \mathbb{R}$,

$$(T_n u)(x) \geq \tilde{Q}(u(x + B)) - \eta_1 \tag{2.10}$$

and

$$(T_n u)(x) \leq \tilde{Q}(u(x - B)) + \eta_1. \tag{2.11}$$

In the special case where $T_n$ is of the form (2.3), with $G_n(y, \cdot) = G(\cdot)$ and $Q_n = Q$, and $Q$ is strictly concave, Assumption 2.4 automatically holds with $\tilde{Q} = Q$. (One can see this by truncating $G$ off a large enough interval which depends on $\eta_1$.)

Our main result asserts the tightness of the distribution functions given by the recursions in (2.1) and (2.3).

THEOREM 2.5. *Let Assumptions 2.1, 2.2 and 2.4 hold, and assume that $F_0(x) = \mathbf{1}_{\{x \geq 0\}}$. Then, the sequence of distribution functions $\{F_n^s\}_{n \geq 0}$ is tight.*

REMARK 2.6. The assumption on $F_0$ in Theorem 2.5 can be relaxed to the assumption that $L(\overline{F}_0) < \infty$ for the function $L$ in Theorem 2.7 below.

The proof of Theorem 2.5 employs a Lyapunov function $L : \mathcal{D} \to \mathbb{R}$ that we introduce next. Choose $\delta_0 > 0$ so that Assumption 2.1 is satisfied, $M$ as in (G2) of Assumption 2.2, $\varepsilon_1 > 0$, and $b > 1$. For $u \in \mathcal{D}$, introduce the function

$$(2.12) \qquad L(u) = \sup_{\{x : u(x) \in (0, \delta_0]\}} \ell(u; x),$$

where

$$(2.13) \qquad \ell(u; x) = \log\left(\frac{1}{u(x)}\right) + \log_b\left(1 + \varepsilon_1 - \frac{u(x - M)}{u(x)}\right)_+.$$

Here, we let $\log 0 = -\infty$ and $(x)_+ = x \vee 0$. If the set on the right-hand side of (2.12) is empty [as it is for $u(x) = \mathbf{1}_{\{x < 0\}}$], we let $L(u) = -\infty$. We will prove the following result in Section 3.

THEOREM 2.7. *Let Assumptions 2.1 and 2.2 hold. There is a choice of parameters $\delta_0, \varepsilon_1, M > 0$ and $b > 1$, such that, if $L(\overline{F}_0) < \infty$, then*

$$(2.14) \qquad \sup_n L(\overline{F}_n) < \infty.$$

Theorem 2.7 implies that, with the given choice of parameters, if $L(\overline{F}_0) < \infty$ and Assumptions 2.1 and 2.2 hold, then, for all $n$ and $x$ with $0 < \overline{F}_n(x) \leq \delta_0$,

$$\log\left(1 + \varepsilon_1 - \frac{\overline{F}_n(x - M)}{\overline{F}_n(x)}\right)_+ \leq (\log b)(C_0 + \log \overline{F}_n(x))$$

for $C_0 = \sup_{n \geq 0} L(\overline{F}_n) < \infty$. In particular, by taking $\delta_1 > 0$ small enough such that the right-hand side in the last inequality is sufficiently negative when $0 < \overline{F}_n(x) \leq \delta_1$, we have the following corollary.

COROLLARY 2.8. *Let Assumptions 2.1 and 2.2 hold. Then, there exists $\delta_1 = \delta_1(C_0, \delta_0, \varepsilon_1, b, M) > 0$ such that, for all $n$,*

$$(2.15) \qquad \overline{F}_n(x) \leq \delta_1 \quad \text{implies} \quad \overline{F}_n(x - M) \geq \left(1 + \frac{\varepsilon_1}{2}\right)\overline{F}_n(x).$$

The inequality (2.15) is sufficient to imply the tightness of $\{F_n^s\}_{n \geq 0}$, irrespective of the precise form of the operator $T_n$. This is shown in the following proposition.

PROPOSITION 2.9. *Suppose that (2.15) holds for all $n$ and some choice of $\delta_1, M, \varepsilon_1 > 0$. Also, suppose that each $T_n$ satisfies Assumption 2.4, and that $\overline{F}_{n+1} = T_n \overline{F}_n$. Then, the sequence of distributions $\{F_n^s\}_{n \geq 0}$ is tight.*

Theorem 2.5 follows directly from Corollary 2.8 and Proposition 2.9. Since the proof of Theorem 2.7, and hence of Corollary 2.8, is considerably longer than that of Proposition 2.9, we prove Proposition 2.9 here and postpone the proof of Theorem 2.7 until Section 3.

The main step in showing Proposition 2.9 is given by Lemma 2.10. It says, in essence, that if $\overline{F}_n$ is "relatively flat" somewhere away from 0 or 1, then $\overline{F}_{n-1}$ is "almost as flat" at a nearby location, where its value is also smaller by a fixed factor $\gamma < 1$. The proof of the lemma will be postponed until after we complete the proof of Proposition 2.9.

LEMMA 2.10. *Suppose that (2.15) holds for all $n$ under some choice of $\delta_1, M, \varepsilon_1 > 0$. Also, suppose that each $T_n$ satisfies Assumption 2.4, and that $\overline{F}_{n+1} = T_n \overline{F}_n$. For fixed $\eta_0 \in (0, 1)$, there exist a constant $\gamma = \gamma(\eta_0) < 1$ and a continuous function $f(t) = f_{\eta_0}(t) : [0, 1] \to [0, 1]$, with $f(t) \to_{t \to 0} 0$ such that for any $\varepsilon \in (0, (1 - \eta_0)/\eta_0)$, $\eta \in [\delta_1, \eta_0]$, and large enough $N_1 = N_1(\varepsilon)$, the following holds. If $M' \geq M$ and, for given $n$ and $x$, $\overline{F}_n(x) > \delta_1$,*

$$(2.16) \qquad \overline{F}_n(x - M') \leq (1 + \varepsilon)\overline{F}_n(x)$$

*and*

$$(2.17) \qquad \overline{F}_n(x - M') \leq \eta,$$

*then*

$$(2.18) \qquad \overline{F}_{n-1}(x + N_1 - M') \leq (1 + f(\varepsilon))\overline{F}_{n-1}(x - N_1)$$

*and*

$$(2.19) \qquad \overline{F}_{n-1}(x + N_1 - M') \leq \gamma \eta.$$

PROOF OF PROPOSITION 2.9 ASSUMING LEMMA 2.10. Fix an $\eta_0 \in (0,1)$. We will show the existence of an $\hat{\varepsilon}_0 = \hat{\varepsilon}_0(\eta_0) > 0$, an $n_0$ and an $\hat{M}$, such that if $n > n_0$ and $\overline{F}_n(x - \hat{M}) \leq \eta_0$, then

$$\overline{F}_n(x - \hat{M}) \geq (1 + \hat{\varepsilon}_0)\overline{F}_n(x). \tag{2.20}$$

This implies the claimed tightness in the statement of the proposition.

The proof of (2.20) is by contradiction, and consists of repeatedly applying Lemma 2.10 until a small enough value of $\overline{F}_n$, where the function is "relatively flat," is achieved, which will contradict (2.15). The presence of the uniform bound $\gamma$ in (2.19) will play an important role in the computations. We proceed by defining $\sigma_i$, $i \geq 0$, by $\sigma_0 = \eta_0$ and $\sigma_i = \gamma\sigma_{i-1}$. Since $\gamma < 1$, the sequence $\sigma_i$ decreases to 0. Set

$$n_0 = \min\{i : \sigma_i < \delta_1\}$$

and specify a sequence $\hat{\varepsilon}_i > 0$, $i \geq 0$, so that $\hat{\varepsilon}_{n_0} < \varepsilon_1/4$ and $\hat{\varepsilon}_i = f(\hat{\varepsilon}_{i-1}) \leq \varepsilon_1/4$. [This is always possible by our assumption that $f(t) \to_{t \to 0} 0$.] Also, set $M_{n_0} = M$ and, for $i \in \{1, \ldots, n_0\}$, set

$$M_{n_0-i} = M_{n_0-i+1} + 2N_1(\hat{\varepsilon}_{n_0-i+1})$$

and $\hat{M} = M_0$.

Suppose now that (2.20) does not hold for some $x$ and $n > n_0$, with $\overline{F}_n(x - \hat{M}) \leq \eta_0 = \sigma_0$. Then,

$$\overline{F}_n(x - \hat{M}) \leq (1 + \hat{\varepsilon}_0)\overline{F}_n(x) \wedge \sigma_0,$$

with $\overline{F}_n(x) > \delta_1$. [The last inequality follows automatically from $\hat{M} \geq M$ and (2.15).] In particular, (2.16) and (2.17) hold with $M' = \hat{M}$, $\varepsilon = \hat{\varepsilon}_0$ and $\eta = \sigma_0$. Applying Lemma 2.10, one concludes that

$$\overline{F}_{n-1}(x + N_1(\hat{\varepsilon}_0) - \hat{M}) \leq (1 + \hat{\varepsilon}_1)\overline{F}_{n-1}(x - N_1(\hat{\varepsilon}_0))$$

and

$$\overline{F}_{n-1}(x + N_1(\hat{\varepsilon}_0) - \hat{M}) \leq \gamma\sigma_0 = \sigma_1.$$

Setting $y = x - N_1(\hat{\varepsilon}_0)$, it follows that there exists a point $y$, such that

$$\overline{F}_{n-1}(y - M_1) \leq (1 + \hat{\varepsilon}_1)\overline{F}_{n-1}(y) \wedge \sigma_1,$$

where $M_1 = \hat{M} - 2N_1(\hat{\varepsilon}_0) \geq M$ by construction.

When $\overline{F}_{n-1}(y) \leq \delta_1$, this contradicts (2.15) because $\hat{\varepsilon}_1 < \varepsilon_1/2$. When $\overline{F}_{n-1}(y) > \delta_1$, repeat the above procedure $n_1$ times (with $n_1 \leq n_0$) to show that there exists a point $y'$ such that

$$\overline{F}_{n-n_1}(y' - M_{n_1}) \leq \delta_1, \qquad \overline{F}_{n-n_1}(y' - M_{n_1}) \leq (1 + \hat{\varepsilon}_{n_1})\overline{F}_{n-n_1}(y').$$

This contradicts (2.15), because $\overline{F}_{n-n_1}(y') \leq \overline{F}_{n-n_1}(y' - M_{n_1})$, $M_{n_1} \geq M$ and $\hat{\varepsilon}_{n_1} \leq \varepsilon_1/4$. $\square$

We now prove Lemma 2.10. The argument for (2.18) consists of two main steps, where one first shows the inequality (2.24) below, and then uses this to show (2.18). The inequality (2.24) follows with the aid of the properties in Assumption 2.4, which allow us to approximate the operator $T_n$ by the pointwise transformation $\tilde{Q}$, after an appropriate translation. The inequality (2.9) is then employed to absorb the coefficient $(1 + 2\varepsilon)$ in (2.24) into the argument of $\tilde{Q}$, from which (2.18) will follow after inverting $\tilde{Q}$. The argument for (2.19) also uses one direction of Assumption 2.4 to bound $T_n$ by $\tilde{Q}$; one then inverts $\tilde{Q}$ to obtain (2.19).

PROOF OF LEMMA 2.10. We first demonstrate (2.18). Suppose (2.16) and (2.17) hold for some $x$ with $\overline{F}_n(x) > \delta_1$; one then also has $\overline{F}_n(x - M') > \delta_1$. Let $\tilde{Q}$ be as in Assumption 2.4. By (2.10), (2.11) and (2.17), for any $\eta_1 > 0$, there exists $B = B(\eta_1)$, such that

$$(2.21) \qquad \overline{F}_n(x - M') \geq \tilde{Q}(\overline{F}_{n-1}(x + B - M')) - \eta_1$$

and

$$(2.22) \qquad \overline{F}_n(x) \leq \tilde{Q}(\overline{F}_{n-1}(x - B)) + \eta_1.$$

By (2.22), since $\overline{F}_n(x) > \delta_1$,

$$(2.23) \qquad \tilde{Q}(\overline{F}_{n-1}(x - B)) \geq \left(1 - \frac{\eta_1}{\delta_1}\right) \overline{F}_n(x).$$

On the other hand, by (2.21) and (2.16),

$$(1+\varepsilon)\overline{F}_n(x) \geq \tilde{Q}(\overline{F}_{n-1}(x+B-M')) - \eta_1 > \tilde{Q}(\overline{F}_{n-1}(x+B-M')) - \frac{\eta_1 \overline{F}_n(x)}{\delta_1}.$$

Combining this with (2.23), it follows that, for $\eta_1 < \delta_1$,

$$(1 + \varepsilon + c(\varepsilon, \eta_1, \delta_1))\tilde{Q}(\overline{F}_{n-1}(x - B)) \geq \tilde{Q}(\overline{F}_{n-1}(x + B - M')),$$

where

$$c(\varepsilon, \eta_1, \delta_1) = \frac{1 + \varepsilon + \eta_1/\delta_1}{1 - \eta_1/\delta_1} - 1 - \varepsilon$$

and, in particular, $c(\varepsilon, \eta_1, \delta_1) \to_{\eta_1 \to 0} 0$. Therefore, for any $\eta_1$ with $c(\varepsilon, \eta_1, \delta_1) < \varepsilon$,

$$(2.24) \qquad (1+2\varepsilon)\tilde{Q}(\overline{F}_{n-1}(x - B)) \geq \tilde{Q}(\overline{F}_{n-1}(x + B - M')).$$

We now choose $N_1 = B$. To prove (2.18), we can assume that $x - B > x + B - M'$. [Otherwise, (2.18) is trivial.] Since $\overline{F}_n(x) > \delta_1$, it follows from (2.22) that, if $\eta_1 < \delta_1/2$, then

$$(2.25) \qquad \tilde{Q}(\overline{F}_{n-1}(x + B - M')) \geq \tilde{Q}(\overline{F}_{n-1}(x - B)) \geq \delta_1/2 > 0.$$

On the other hand, from (2.17) and (2.21),

$$(2.26) \qquad \eta \geq \overline{F}_n(x - M') \geq \tilde{Q}(\overline{F}_{n-1}(x + B - M')) - \eta_1.$$

In particular, for each $\eta_1 < (1 - \eta)/2$ and $\delta' = \delta'(\eta)$ chosen so that

$$\delta' = \tilde{Q}^{-1}\left(\frac{\delta_1}{2}\right) \wedge \frac{1 - \eta}{2} > 0,$$

$\overline{F}_{n-1}(x + B - M') \geq \delta'$ and $\tilde{Q}(\overline{F}_{n-1}(x + B - M')) \leq 1 - \delta'$. Applying (2.9) together with (2.24), one concludes that

$$(1 + f(\varepsilon))\overline{F}_{n-1}(x - B) \geq \overline{F}_{n-1}(x + B - M'),$$

with the function $f(\varepsilon) := g_{\delta'}(2\varepsilon) \to_{\varepsilon \to 0} 0$. The inequality (2.18) follows since $N_1 = B$.

Note that by (2.26), for any $\eta_1 > 0$,

$$\overline{F}_{n-1}(x + B - M') \leq \tilde{Q}^{-1}(\overline{F}_n(x - M') + \eta_1) \leq \tilde{Q}^{-1}(\eta + \eta_1).$$

The inequality (2.19) follows from this and property (T1), by choosing $\eta_1$ small enough so that $\gamma = \sup_{\eta \in [\delta_1, \eta_0]} \tilde{Q}^{-1}(\eta + \eta_1)/\eta < 1$. □

**3. Proof of Theorem 2.7.** This section is devoted to proving Theorem 2.7. In order to prove the result, we will show the following minor variation.

THEOREM 3.1. *Let Assumptions 2.1 and 2.2 hold. There is a choice of parameters $\delta_0, \varepsilon_1, M, C_1 > 0$ and $b > 1$, with the property that if $L(\overline{F}_{n+1}) \geq C$ for some $n$ and some $C > C_1$, then $L(\overline{F}_n) \geq C$.*

PROOF OF THEOREM 2.7 ASSUMING THEOREM 3.1. If $\sup_n L(\overline{F}_n) = \infty$, then for any $C$, one can choose $n$ such that $L(\overline{F}_n) \geq C$. For $C > C_1$, it follows by Theorem 3.1 that $L(\overline{F}_0) \geq C$. Since $C$ can be made arbitrarily large, one must have $L(\overline{F}_0) = \infty$, which contradicts the assumption that $L(\overline{F}_0) < \infty$. □

In what follows, for a function $f$ on $\mathbb{R}$ and any $x \in \mathbb{R}$, we set $f(x)^- = \lim_{y \nearrow x} f(y)$, when this limit exists. We define, for $\varepsilon > 0$, $x_1 \in \mathbb{R}$, $x_2 = x_1 - M$, and $u \in \mathcal{D}$,

$$q_1 = q_1(u, \varepsilon, x_1) = \inf\{y > 0 : u(x_2 - y) \geq (1 + 8\varepsilon)u(x_1 - y)\}$$

and

$$r_1 = r_1(u, \varepsilon, x_1)$$
$$= \begin{cases} q_1, & \text{if } u(x_2 - q_1)^- \geq u(x_1 - q_1 - M/2)/(1 - 4\varepsilon), \\ q_1 - \dfrac{M}{2}, & \text{otherwise.} \end{cases}$$

Possibly, $q_1(u, \varepsilon, x_1) = \infty$, in which case we set $r_1(u, \varepsilon, x_1) = \infty$. Intuitively, $x_1 - q_1$ is the first point to the left of $x_1$ where $u$ is "very nonflat." [We are interpreting $u$ to be "very nonflat" at $x$ if the ratio $u(x - M)/u(x)$ is not close to 1.] We have defined $r_1$ so that $u$ is "very nonflat" at all points in $[x_1 - r_1 - M/2, x_1 - r_1)$; more detail will be given in the proof of Lemma 3.5.

The proof of Theorem 3.1 is based on the following two propositions. The first allows one to "deconvolve" the ⊛ operation and maintain a certain amount of "nonflatness." In the remainder of the section we will implicitly assume that $aM > 100$ and that $M > 2M_0$, where $M_0$ is as in Assumption (G2).

PROPOSITION 3.2. *Let Assumption 2.2 hold. For a given $u \in \mathcal{D}$, $n$, $x_1 \in \mathbb{R}$ and $\varepsilon' \in (0, 1/64)$, and for $x_2 = x_1 - M$, assume that*

$$(3.1) \qquad (\overline{G}_n \circledast u)(x_2) < (1 + \varepsilon')(\overline{G}_n \circledast u)(x_1).$$

*Then, at least one of the following two statements holds for each $\delta > 0$:*

$(3.2) \quad u(x_2 - y) \leq (1 + \varepsilon' + \delta)u(x_1 - y), \qquad \text{some } y \leq M \wedge r_1(u, \varepsilon', x_1)$

$(3.3) \quad u(x_2 - y) \leq (1 + \varepsilon' - \delta e^{ay/4})u(x_1 - y), \qquad \text{some } y \in (M, r_1(u, \varepsilon', x_1)].$

The second proposition controls the Lyapunov function $\ell$ around "nonflat" points.

PROPOSITION 3.3. *Let Assumption 2.2 hold. For a given $u \in \mathcal{D}$, $n$, $x_1 \in \mathbb{R}$ and $\varepsilon \in [0, \varepsilon_1)$, with $\varepsilon_1 \leq 1/64$, and for $x_2 = x_1 - M$, assume that*

$$(3.4) \qquad (\overline{G}_n \circledast u)(x_2) = (1 + \varepsilon)(\overline{G}_n \circledast u)(x_1).$$

*Choose $\delta < \kappa(\varepsilon_1 - \varepsilon)$ and $\varepsilon' = \varepsilon + \delta/2 < \varepsilon_1$, where $\kappa \in (0, 1)$. Then, for small enough $\kappa > 0$ and for $b > 1$, neither depending on $u, x_1, x_2, \delta$ or $\varepsilon$, the following hold:*

(a) *If (3.2) is satisfied, then there exists $x_1' \geq x_1 - M$ such that*

$$(3.5) \qquad \ell(u; x_1') - \ell(\overline{G}_n \circledast u; x_1) \geq -2\left(\varepsilon_1 + \varepsilon_0 + \dfrac{\delta}{(\varepsilon_1 - \varepsilon)\log b}\right).$$

(b) *If (3.3) is satisfied, then there exists $x_1' \leq x_1 - M$ such that*

$$\ell(u; x_1') - \ell(\overline{G}_n \circledast u; x_1)$$
(3.6)
$$\geq -6(\varepsilon_1(x_1 - x_1') + \varepsilon_0) + \frac{a(x_1 - x_1') + 4\log(\delta/(\varepsilon_1 - \varepsilon))}{4 \log b}.$$

We will also employ the following lemma, which allows us to avoid checking the condition $x_2 \leq 2x_1$ in (2.5).

LEMMA 3.4. *Let Assumption 2.1 hold. Suppose for small enough $\delta_0 > 0$ that $0 < x_1 < x_2$ satisfy*

(3.7) $$Q_n(x_1) \leq \delta_0 \quad \text{and} \quad Q_n(x_2) \leq \tfrac{3}{2} Q_n(x_1).$$

*Then, for $c^*$ and $\theta^*$ chosen as in the assumption,*

(3.8) $$\frac{x_2}{x_1} \leq \frac{Q_n(x_2)}{Q_n(x_1)} [1 + c^*(Q_n(x_1))^{\theta^*}].$$

PROOF OF LEMMA 3.4. Since $Q_n$ is increasing, $x_1 \leq \delta_0$ by (2.4) and (3.7); the main inequality in (2.5) will therefore hold if $x_2 \wedge 2x_1$ is substituted there for $x_2$. Together with (3.7), this implies that

$$\frac{x_2 \wedge 2x_1}{x_1} \leq \frac{3}{2}[1 + c^*(Q_n(x_1))^{\theta^*}],$$

which is $< 2$ for small enough $\delta_0$ and, hence, $x_2 < 2x_1$. One can now employ (2.5) to obtain (3.8). $\square$

The proof of Proposition 3.2 requires a fair amount of work. We therefore first demonstrate Theorem 3.1 assuming Propositions 3.2 and 3.3, and then afterward demonstrate both propositions.

PROOF OF THEOREM 3.1. Assume that the constant $\varepsilon_1$ in (2.13) satisfies $\varepsilon_1 < (\log m)/100$, where $m \in (1, 2]$ is as in (2.4). We will show that, for $C$ large enough and $L(\overline{F}_{n+1}) \geq C$,

(3.9) $$L(Q_n(\overline{F}_n)) - L(\overline{G}_n \circledast Q_n(\overline{F}_n)) \geq -\frac{\log m}{4}$$

and

(3.10) $$L(\overline{F}_n) - L(Q_n(\overline{F}_n)) \geq \frac{\log m}{2}.$$

Since $L(\overline{F}_{n+1}) = L(\overline{G}_n \circledast Q_n(\overline{F}_n))$, it will follow from (3.9) and (3.10) that $L(\overline{F}_n) \geq L(\overline{F}_{n+1})$, which implies Theorem 3.1.

Recall the constants $\theta^*$ and $a$ from Assumptions 2.1 and 2.2, and the constant $\kappa$ from Proposition 3.3. We fix $\delta_0 > 0$ small enough so that the conclusions of Proposition 3.3 and Lemma 3.4 hold. We also choose $b$ close enough to 1, with $e^{\theta^*} > b > 1$, and $M > 100/a$ large enough, such that the conclusion of Proposition 3.3 and all the following conditions hold:

$$\text{(3.11)} \qquad \frac{(\log b)(\log m)}{20} \leq \kappa,$$

$$\text{(3.12)} \qquad \frac{a}{8 \log b} > 6\varepsilon_1 + \frac{6 \log m}{100M},$$

$$\text{(3.13)} \qquad \frac{aM}{8} > -\log\left(\frac{(\log b)(\log m)}{20}\right).$$

We begin with the proof of (3.9). Propositions 3.2 and 3.3 provide the main ingredients. Choose $x_1$ such that

$$(\overline{G}_n \circledast Q_n(\overline{F}_n))(x_1) \leq \delta_0, \qquad \ell(\overline{G}_n \circledast Q_n(\overline{F}_n); x_1) > C - 1$$

and

$$\text{(3.14)} \qquad \ell(\overline{G}_n \circledast Q_n(\overline{F}_n); x_1) > L(\overline{G}_n \circledast Q_n(\overline{F}_n)) - (\log m)/10.$$

This is always possible since $L(\overline{F}_{n+1}) = L(\overline{G}_n \circledast Q_n(\overline{F}_n)) \geq C$. Choose $\varepsilon \in [0, \varepsilon_1)$ such that (3.4) in Proposition 3.3 holds for $u = Q_n(\overline{F}_n)$. Also, set $\delta = (\log b)(\log m)(\varepsilon_1 - \varepsilon)/40$ and note that, due to (3.11), $\delta < \kappa(\varepsilon_1 - \varepsilon)$. Applying Proposition 3.2, with $Q_n(\overline{F}_n)$ playing the role of $u$ there and with $\varepsilon' = \varepsilon + \delta/2$, either (3.2) or (3.3) must hold.

Suppose that (3.2) holds, and set $\alpha_1 = 2(\varepsilon_1 + \varepsilon_0) + (\log m)/20 \leq (\log m)/10$, where $\varepsilon_0$ is given in (2.7). Then, by (3.5) of Proposition 3.3, there exists $x_1' \geq x_1 - M$ such that

$$\text{(3.15)} \qquad \ell(Q_n(\overline{F}_n); x_1') - \ell(\overline{G}_n \circledast Q_n(\overline{F}_n); x_1) \geq -\alpha_1.$$

In particular, $\ell(Q_n(\overline{F}_n); x_1') \geq C - 1 - \alpha_1$ and, hence, by the definition of $\ell$,

$$\text{(3.16)} \qquad -\log Q_n(\overline{F}_n(x_1')) \geq C - 1 - \alpha_1 - \log_b(1 + \varepsilon_1)$$

with the right-hand side being greater than $-\log \delta_0$ if $C$ is large enough. Together with (3.15), (3.14) and the definition of $L$, this yields (3.9) when (3.2) holds.

Suppose now that (3.3) holds. Then, again by Proposition 3.3, (3.6) implies that

$$\text{(3.17)} \quad \begin{aligned} &\ell(Q_n(\overline{F}_n); x_1') - \ell(\overline{G}_n \circledast Q_n(\overline{F}_n); x_1) \\ &\qquad \geq -6(\varepsilon_1(x_1 - x_1') + \varepsilon_0) + \frac{a(x_1 - x_1')}{4 \log b} + \frac{\log((\log b)(\log m)/20)}{\log b}. \end{aligned}$$

Since $\varepsilon_0 = (\log m)/100$, it follows by (3.12) and (3.13) that the right-hand side of (3.17) is nonnegative. Further, by exactly the same argument as in the case that (3.2) holds (after replacing $-\alpha_1$ by 0), one deduces from (3.17) that

$$-\log Q_n(\overline{F}_n(x_1')) \geq C - 1 - \log_b(1 + \varepsilon_1) > -\log \delta_0.$$

Together with (3.14), (3.17) and the definition of $L$, this yields (3.9) when (3.3) holds and, hence, completes the proof of (3.9).

The proof of (3.10) is based on the properties of $Q_n$ given in Assumption 2.1. The basic idea is that when $\overline{F}_n(x_1)$ is sufficiently small, $\ell(\overline{F}_n; x_1) - \ell(Q_n(\overline{F}_n); x_1)$ will be almost $\log m$ because: (a) the difference of the first components of $\ell$ contributes

$$\log(Q_n(\overline{F}_n(x_1))/\overline{F}_n(x_1)) \geq \log m$$

on account of (2.4) and (b) the difference of the second components of $\ell$ is negligible, since $\overline{F}_n(x_2)/\overline{F}_n(x_1)$ is not much larger than $Q_n(\overline{F}_n(x_2))/Q_n(\overline{F}_n(x_1))$, on account of (2.5).

To justify this reasoning, first note that, by (3.9), we already know that $L(Q_n(\overline{F}_n)) > C - (\log m)/3 > 0$. So, there exists an $x_1$ with $Q_n(\overline{F}_n(x_1)) \leq \delta_0$ and

$$(3.18) \quad \ell(Q_n(\overline{F}_n); x_1) > \max\left(C - \frac{\log m}{3}, L(Q_n(\overline{F}_n)) - \frac{\log m}{10}\right).$$

Since the right-hand side is $> -\infty$, one has, for $x_2 = x_1 - M$,

$$(3.19) \qquad Q_n(\overline{F}_n(x_2)) \leq (1 + \varepsilon_1) Q_n(\overline{F}_n(x_1)).$$

So, (3.7) is satisfied, with $\overline{F}_n(x_i)$ in place of $x_i$, since $\varepsilon_1 < 1/2$. Consequently, by Lemma 3.4,

$$(3.20) \qquad \frac{\overline{F}_n(x_2)}{\overline{F}_n(x_1)} \leq \frac{Q_n(\overline{F}_n(x_2))}{Q_n(\overline{F}_n(x_1))}[1 + c^*(Q_n(\overline{F}_n(x_1)))^{\theta^*}].$$

Now, set $Q_n(\overline{F}_n(x_1)) = q$ and $Q_n(\overline{F}_n(x_2))/Q_n(\overline{F}_n(x_1)) = 1 + \varepsilon$, which is less than $1 + \varepsilon_1$. By (2.4), $\overline{F}_n(x_1) \leq \delta_0$ holds and, hence, $Q_n(\overline{F}_n(x_1)) \geq m\overline{F}_n(x_1)$. So,

$$\ell(\overline{F}_n; x_1) = -\log \overline{F}_n(x_1) + \log_b\left(1 + \varepsilon_1 - \frac{\overline{F}_n(x_2)}{\overline{F}_n(x_1)}\right)_+$$

$$\geq -\log Q_n(\overline{F}_n(x_1)) + \log m + \log_b(\varepsilon_1 - \varepsilon - 2c^* q^{\theta^*})_+.$$

Because $\ell(Q_n(\overline{F}_n); x_1) > C - (\log m)/3$,

$$q \leq e^{-(C-(\log m)/3)}(\varepsilon_1 - \varepsilon)^{1/\log b}.$$

Since $\theta^*/\log b > 1$, it follows from the previous two displays and the definition of $\ell$ that

$$\ell(\overline{F}_n; x_1) \geq -\log Q_n(\overline{F}_n(x_1)) + \log m$$
$$+ \log_b(\varepsilon_1 - \varepsilon - 2c^* e^{-\theta^*(C-(\log m)/3)}(\varepsilon_1 - \varepsilon)^{\theta^*/\log b})_+$$
$$\geq \ell(Q_n(\overline{F}_n); x_1) + \log m + \log_b(1 - 2c^* e^{-\theta^*(C-(\log m)/3)})_+.$$

Choosing $C$ large enough such that

$$\log_b(1 - 2c^* e^{-\theta^*(C-(\log m)/3)})_+ \geq -\frac{\log m}{10}$$

and using (3.18), we get

$$L(\overline{F}_n) \geq \ell(\overline{F}_n; x_1) \geq \ell(Q_n(\overline{F}_n); x_1) + \tfrac{9}{10}\log m$$
$$\geq L(Q_n(\overline{F}_n)) + \tfrac{1}{2}\log m.$$

This implies (3.10). □

We now prove Proposition 3.2, which was used in the proof of Theorem 3.1. Much of the work is contained in the following lemma, whose demonstration we postpone until after that of the proposition.

LEMMA 3.5. *Suppose $\overline{G}_n$ satisfy Assumption 2.2, and for given $\varepsilon' > 0$, $u \in \mathcal{D}$, $x_1 \in \mathbb{R}$ and $x_2 = x_1 - M$, that (3.1) holds. Then,*

$$\int_{x_2-r_1}^{\infty} u(y)\, d\overline{G}_n(y+M, x_2-y) = \int_{x_1-r_1}^{\infty} u(y-M)\, d\overline{G}_n(y, x_1-y)$$
(3.21)
$$< (1+\varepsilon') \int_{x_1-r_1}^{\infty} u(y)\, d\overline{G}_n(y, x_1-y),$$

*where $r_1 = r_1(u, \varepsilon', x_1)$.*

The integration is to be interpreted as being over the parameter $y$. Since $\overline{G}_n(y, x-y)$ is increasing in $y$ for each fixed $x$, the integrals are well defined as Lebesgue–Stieltjes integrals. Here and later on, we use the convention

$$\int_a^b f(y)\, dg(y) = \int_{-\infty}^{\infty} \mathbf{1}_{\{y \in [a,b)\}} f(y)\, dg(y).$$

PROOF OF PROPOSITION 3.2. In the proof we will omit the index $n$ from $\overline{G}_n$, writing $\overline{G}$ since the estimates that are used do not depend on $n$.

Suppose that $r_1 = r_1(u, \varepsilon', x_1) \leq M$. If (3.2) does not hold for a given $\delta > 0$, then

$$u(y-M) \geq (1+\varepsilon'+\delta) u(y) \quad \text{for all } y \geq x_1 - r_1,$$

which contradicts (3.21). So, to prove Proposition 3.2, it remains only to consider the case where $r_1 > M$.

Assume now that for a given $\delta > 0$, neither (3.2) nor (3.3) holds. We will show that this again contradicts (3.21). Decompose the left side of (3.21) into the integrals over $[x_2 - r_1, x_2 - M)$ and $[x_2 - M, \infty)$, which we denote by $A_1$ and $A_2$. Since (3.2) is violated,

$$
\begin{aligned}
A_2 &= (1 + \varepsilon' + \delta) \int_{x_2 - M}^{\infty} u(y) \, d\overline{G}(y + M, x_2 - y) \\
&\geq (1 + \varepsilon' + \delta) \int_{x_1 - M}^{\infty} u(y) \, d\overline{G}(y, x_1 - y).
\end{aligned}
$$
(3.22)

Since (3.3) is violated,

$$
\begin{aligned}
A_1 &= \int_{x_1 - r_1}^{x_1 - M} u(y - M) \, d\overline{G}(y, x_1 - y) \\
&\geq (1 + \varepsilon') \int_{x_1 - r_1}^{x_1 - M} u(y) \, d\overline{G}(y, x_1 - y) \\
&\quad - \delta \int_{x_1 - r_1}^{x_1 - M} u(y) e^{a(x_1 - y)/4} \, d\overline{G}(y, x_1 - y).
\end{aligned}
$$

Hence,

$$
\begin{aligned}
A_1 + A_2 &- (1 + \varepsilon') \int_{x_1 - r_1}^{\infty} u(y) \, d\overline{G}(y, x_1 - y) \\
&\geq \delta \bigg[ \int_{x_1 - M}^{\infty} u(y) \, d\overline{G}(y, x_1 - y) \\
&\quad - \int_{x_1 - r_1}^{x_1 - M} u(y) e^{a(x_1 - y)/4} \, d\overline{G}(y, x_1 - y) \bigg].
\end{aligned}
$$
(3.23)

We will show that the right-hand side of (3.23) is nonnegative, which will contradict (3.21) and demonstrate the proposition.

We will bound the second integral on the right-hand side of (3.23). The basic idea will be to control the growth of $u(y)$ as $y$ decreases; an exponential bound on this follows from the definition of $r_1$. The exponential bound in Assumption (G2) on the tail of $G_n$ dominates this rate, and allows us to bound the corresponding integral by the first integral on the right-hand side of (3.23).

Proceeding with the argument for this, one can check that

$$
\begin{aligned}
\int_{x_1 - r_1}^{x_1 - M} &u(y) e^{a(x_1 - y)/4} \, d\overline{G}(y, x_1 - y) \\
&= \sum_{k=1}^{\infty} \int_{x_1 - r_1}^{x_1 - M} \mathbf{1}_{\{y - x_1 \in (-(k+1)M, -kM]\}} u(y) e^{a(x_1 - y)/4} \, d\overline{G}(y, x_1 - y)
\end{aligned}
$$

(3.24)
$$\leq \sum_{k=1}^{\infty} e^{a(k+1)M/4}$$
$$\cdot \int_{x_1-M}^{x_1} \mathbf{1}_{\{y-x_1>kM-r_1\}} u(y-kM) \, d\overline{G}(y-kM, x_1-y+kM).$$

It follows from the definitions of $q_1$ and $r_1$ that if $y - kM \geq x_1 - r_1$ (and hence, $y - kM \geq x_1 - q_1$),

$$u(y-kM) \leq (1+8\varepsilon')u(y-(k-1)M) \leq \cdots \leq (1+8\varepsilon')^k u(y) \leq e^{akM/8} u(y),$$

where the last inequality holds since $aM > 100$ was assumed. The last expression in (3.24) is therefore at most

$$\sum_{k=1}^{\infty} e^{a(k+1)M/8} e^{akM/4} \int_{x_1-M}^{x_1} u(y) \, d\overline{G}(y-kM, x_1-y+kM).$$

Since $u(y)$ is decreasing in $y$ and $u(x_1 - M) \leq e^{aM/8} u(x_1)$, this is at most

$$\sum_{k=1}^{\infty} e^{akM/2} e^{aM/4} u(x_1) \overline{G}(x_1-kM, kM).$$

Applying (G2) $k$ times to $\overline{G}(x_1 - kM, kM)$, summing over $k$, and applying (G2) again gives the upper bound

$$\frac{e^{-aM/4}}{1-e^{-aM/2}} u(x_1) \overline{G}(x_1, 0) \leq \frac{e^{-aM/4}}{(1-e^{-aM/2})^2} u(x_1)(\overline{G}(x_1, 0) - \overline{G}(x_1-M, M)).$$

Since $aM > 100$, this is at most

(3.25) $$\int_{x_1-M}^{x_1} u(y) \, d\overline{G}(y, x_1-y) \leq \int_{x_1-M}^{\infty} u(y) \, d\overline{G}(y, x_1-y).$$

Consequently, by the inequalities (3.24) through (3.25),

(3.26) $$\int_{x_1-r_1}^{x_1-M} u(y) e^{a(x_1-y)/4} \, d\overline{G}(y, x_1-y) \leq \int_{x_1-M}^{\infty} u(y) \, d\overline{G}(y, x_1-y).$$

This shows the right-hand side of (3.23) is nonnegative, and completes the proof of Proposition 3.2. □

PROOF OF LEMMA 3.5. In the proof we will omit the index $n$ from $\overline{G}_n$, writing $\overline{G}$ since the estimates that are used do not depend on $n$.

Since $x_2 = x_1 - M$, the equality in (3.21) is immediate. In order to demonstrate the inequality in (3.21), we note that

$$(1+\varepsilon')(\overline{G} \circledast u)(x_1) > (\overline{G} \circledast u)(x_2)$$

$$= -\int_{-\infty}^{\infty} \overline{G}(y, x_2 - y) \, du(y)$$

(3.27)
$$\geq -\int_{-\infty}^{\infty} \overline{G}(y + M, x_2 - y) \, du(y)$$

$$= +\int_{-\infty}^{\infty} u(y) \, d\overline{G}(y + M, x_2 - y),$$

where the first inequality follows from (3.1) and the second inequality from Assumption (G1). Subtracting the right and left-hand sides of (3.21) from the first and last terms in (3.27), it therefore suffices to show that

(3.28) $$\int_{-\infty}^{x_2-r_1} u(y) \, d\overline{G}(y + M, x_2 - y) \geq (1 + \varepsilon') \int_{-\infty}^{x_1-r_1} u(y) \, d\overline{G}(y, x_1 - y),$$

where $r_1 = r_1(u, \varepsilon', x_1)$.

We will show (3.28) by partitioning $(-\infty, x_2 - r_1)$ into appropriate subintervals. Starting from $q_1$ and $r_1$ as defined previously, we set, for $k > 1$,

$$q_k = q_k(u, \varepsilon', x_1) = \inf\{y > r_{k-1} + M : u(x_2 - y) \geq (1 + 8\varepsilon') u(x_1 - y)\}$$

and

$$r_k = r_k(u, \varepsilon', x_1) = \begin{cases} q_k, & \text{if } u(x_2 - q_k)^- \geq \dfrac{u(x_1 - q_k - M/2)}{1 - 4\varepsilon}, \\ q_k - \dfrac{M}{2}, & \text{otherwise.} \end{cases}$$

For $K = \sup\{k : q_k < \infty\}$, $q_1, q_2, \ldots, q_K$ are intuitively the successive points at which $u$ is "very nonflat" and which are sufficiently separated; $u$ is "very nonflat" at all points in $[x_1 - r_k - M/2, x_1 - r_k)$.

For $k \leq K$ and $i \in \{1, 2\}$, we set

$$A_k^i = [x_i - r_k - M/2, x_i - r_k), \qquad B_k^i = [x_i - r_{k+1}, x_i - r_k - M/2),$$

and let $q_{K+1} = r_{K+1} = \infty$. By summing over $k \in \{1, \ldots, K\}$, (3.28) will follow once we show that

(3.29)
$$\int_{A_k^2} u(y) \, d\overline{G}(y + M, x_2 - y) - (1 + \varepsilon') \int_{A_k^1} u(y) \, d\overline{G}(y, x_1 - y)$$
$$\geq 2\varepsilon' u(x_2 - r_k)^- \overline{G}(x_1 - r_k, r_k)$$

and

(3.30)
$$\int_{B_k^2} u(y) \, d\overline{G}(y + M, x_2 - y) - (1 + \varepsilon') \int_{B_k^1} u(y) \, d\overline{G}(y, x_1 - y)$$
$$\geq -2\varepsilon' u(x_2 - r_k)^- \overline{G}(x_1 - r_k, r_k).$$

We first show (3.29) for $r_k = q_k$. We do this by deriving a uniform lower bound on $u(y)$ for $y \in A_k^2$, and a uniform upper bound on $u(y)$ for $y \in A_k^1$. Trivially, when $y \in A_k^2$, $u(y) \geq u(x_2 - r_k)$ holds. When $y \in A_k^1$,

$$u(y) \leq u(x_1 - r_k - M/2) = u(x_1 - q_k - M/2)$$
$$\leq (1 - 4\varepsilon')u(x_1 - q_k - M)^- = (1 - 4\varepsilon')u(x_2 - r_k)^-.$$

Also, note that $x_2$ and $A_k^2$ are each obtained by translating $x_1$ and $A_k^1$ by $-M$. Therefore,

$$\int_{A_k^2} u(y)\, d\overline{G}(y + M, x_2 - y) - (1 + \varepsilon') \int_{A_k^1} u(y)\, d\overline{G}(y, x_1 - y)$$
$$\geq \left[ \int_{A_k^2} d\overline{G}(y + M, x_2 - y) - (1 + \varepsilon')(1 - 4\varepsilon') \int_{A_k^1} d\overline{G}(y, x_1 - y) \right]$$
$$\times u(x_2 - r_k)^-$$
$$= [1 - (1 + \varepsilon')(1 - 4\varepsilon')]$$
$$\cdot [\overline{G}(x_1 - r_k, r_k) - \overline{G}(x_1 - r_k - M/2, r_k + M/2)]u(x_2 - r_k)^-.$$

By applying Assumption (G2) and $aM > 100$, it follows that this is

$$\geq [1 - (1 + \varepsilon')(1 - 4\varepsilon')][1 - e^{-aM/2}]\overline{G}(x_1 - r_k, r_k)u(x_2 - r_k)^-$$
$$\geq 2\varepsilon' \overline{G}(x_1 - r_k, r_k)u(x_2 - r_k)^-.$$

This shows (3.29) when $r_k = q_k$.

For $r_k = q_k - M/2$, we employ an analogous argument. When $y \in A_k^1$,

$$u(y) \leq u(x_1 - r_k - M/2) = u(x_1 - q_k)$$
$$\leq (1 + 8\varepsilon')^{-1}u(x_2 - q_k)^- = (1 + 8\varepsilon')^{-1}u(x_2 - r_k - M/2)^-.$$

When $y \in A_k^2$,

$$u(y) \geq u(x_2 - r_k) = u(x_2 - q_k + M/2) = u(x_1 - q_k - M/2)$$
$$\geq (1 - 4\varepsilon')u(x_2 - q_k)^- = (1 - 4\varepsilon')u(x_2 - r_k - M/2)^-.$$

Arguing in the same manner as before, we now get

$$\int_{A_k^2} u(y)\, d\overline{G}(y + M, x_2 - y) - (1 + \varepsilon') \int_{A_k^1} u(y)\, d\overline{G}(y, x_1 - y)$$
$$\geq \left[ 1 - 4\varepsilon' - \frac{1 + \varepsilon'}{1 + 8\varepsilon'} \right][1 - e^{-aM/2}]\overline{G}(x_1 - r_k, r_k)$$
$$\times u(x_2 - r_k - M/2)^-$$
$$\geq 2\varepsilon' \overline{G}(x_1 - r_k, r_k)u(x_2 - r_k)^-.$$

RECURSIONS AND TIGHTNESS 25This completes the proof of (3.29).

We still need to show (3.30). Bound the left side of (3.30) using

$$
\begin{aligned}
\int_{B_k^2} u(y)\, d\overline{G}(y+M, x_2 - y) &- (1+\varepsilon') \int_{B_k^1} u(y)\, d\overline{G}(y, x_1 - y) \\
&= \int_{B_k^1} [u(y - M) - (1+\varepsilon')u(y)]\, d\overline{G}(y, x_1 - y) \\
&\geq - \int_{B_k^1} \varepsilon' u(y)\, d\overline{G}(y, x_1 - y) \\
&= -\varepsilon' \bigg[ \int_{x_1 - r_k - M}^{x_1 - r_k - M/2} u(y)\, d\overline{G}(y, x_1 - y) \\
&\qquad + \int_{x_1 - r_{k+1}}^{x_1 - r_k - M} u(y)\, d\overline{G}(y, x_1 - y) \bigg].
\end{aligned}
\tag{3.31}
$$

Using Assumption (G1), one obtains for the first term on the right that

$$
\begin{aligned}
\int_{x_1 - r_k - M}^{x_1 - r_k - M/2} & u(y)\, d\overline{G}(y, x_1 - y) \\
&\leq u(x_1 - r_k - M) \overline{G}(x_1 - r_k - M/2, r_k + M/2) \\
&\leq u(x_2 - r_k) \overline{G}(x_1 - r_k, r_k).
\end{aligned}
\tag{3.32}
$$

We can also bound the second term on the right side of (3.31). Here, because the interval can be quite long, we divide it up into subintervals of length $M$. Let $L$ denote the smallest integer such that $r_k + ML \geq q_{k+1}$ (possibly, $L = \infty$). We have, using (G1) in the second inequality, the definitions of $q_k$ and $r_k$ in the third, and (G2) together with $r_1 \geq -M/2$ and $M/2 \geq M_0$ in the fourth,

$$
\begin{aligned}
\int_{x_1 - r_{k+1}}^{x_1 - r_k - M} u(y)\, d\overline{G}(y, x_1 - y) &\leq \sum_{\ell=1}^{L-1} \int_{x_1 - r_k - M(\ell+1)}^{x_1 - r_k - M\ell} u(y)\, d\overline{G}(y, x_1 - y) \\
&\leq \sum_{\ell=1}^{L-1} u(x_1 - r_k - M(\ell+1)) \overline{G}(x_1 - r_k - M\ell, r_k + M\ell) \\
&\leq \sum_{\ell=1}^{L-1} (1 + 8\varepsilon')^\ell u(x_2 - r_k)^- \overline{G}(x_1 - r_k - M\ell, r_k + M\ell) \\
&\leq \sum_{\ell=1}^{L-1} (1 + 8\varepsilon')^\ell e^{-aM\ell + aM/2} u(x_2 - r_k)^- \overline{G}(x_1 - r_k - M/2, r_k + M/2).
\end{aligned}
$$

Since $aM > 100$ and $\overline{G}(x_1 - r_k - M/2, r_k + M/2) \leq \overline{G}(x_1 - r_k, r_k)$ it follows that
$$\int_{x_1 - r_{k+1}}^{x_1 - r_k - M} u(y) \, d\overline{G}(y, x_1 - y) \leq u(x_2 - r_k)^{-}\overline{G}(x_1 - r_k, r_k).$$

Substituting this and (3.32) into (3.31) yields (3.30), and completes the proof of Lemma 3.5. □

We now prove Proposition 3.3.

PROOF OF PROPOSITION 3.3. Both parts (a) and (b) derive bounds for the two components of $\ell$ in (2.13). Part (a) employs routine estimates; in part (b), the definition of $r_1$ is employed to control the growth of $u(x)$ as $x$ decreases. We will omit the index $n$ from $\overline{G}_n$, writing $\overline{G}$ since the estimates that are used do not depend on $n$.

We first demonstrate part (a). Assume that (3.2) is satisfied for a given $y = \hat{y}$, and set $x_1' = x_1 - \hat{y}$; since $\hat{y} \leq M$, one has $x_1' \geq x_2$. Then, using (G1) in the next to last inequality and (G3) in the last,

$$\begin{aligned}
(\overline{G} \circledast u)(x_2) &= -\int_{-\infty}^{\infty} \overline{G}(y, x_2 - y) \, du(y) \\
&\geq -\int_{x_2}^{\infty} \overline{G}(y, x_2 - y) \, du(y) \\
&\geq u(x_2) \min_{y \geq x_2} \overline{G}(y, x_2 - y) \\
&\geq u(x_2)\overline{G}(x_2, 0) \geq (1 - \varepsilon_0)u(x_2).
\end{aligned}$$

Hence, using (3.4), one obtains

$$(3.33) \qquad u(x_2) \leq \frac{(\overline{G} \circledast u)(x_2)}{1 - \varepsilon_0} = \frac{1 + \varepsilon}{1 - \varepsilon_0}(\overline{G} \circledast u)(x_1).$$

Since $u(x_1') \leq u(x_2)$, it follows from (3.33) that

$$(3.34) \qquad \log\left(\frac{(\overline{G} \circledast u)(x_1)}{u(x_1')}\right) \geq -\log\left(\frac{1 + \varepsilon}{1 - \varepsilon_0}\right) \geq -2(\varepsilon_1 + \varepsilon_0).$$

On the other hand, setting $x_2' = x_1' - M$, and applying (3.2) and (3.4),

$$\log_b\left(\frac{1 + \varepsilon_1 - u(x_2')/u(x_1')}{1 + \varepsilon_1 - (\overline{G} \circledast u)(x_2)/(\overline{G} \circledast u)(x_1)}\right)$$
$$\geq \log_b\left(\frac{\varepsilon_1 - \varepsilon' - \delta}{\varepsilon_1 - \varepsilon}\right) = \log_b\left(1 - \frac{3\delta}{2(\varepsilon_1 - \varepsilon)}\right),$$

which is larger than $-2\delta/(\varepsilon_1 - \varepsilon) \log b$ for $\kappa$ chosen small enough. Together with the bound in (3.34), this implies (3.5).

We still need to demonstrate part (b). Assume that (3.3) is satisfied for a given $y = \hat{y}$. As before, set $x'_i = x_i - \hat{y}$, where we now have $\hat{y} > M$ and, therefore, $x'_1 < x_1 - M = x_2$. By the same reasoning as in part (a), (3.33) continues to hold. Writing $\hat{y}_M = \lceil \hat{y}/M \rceil - 1 \geq 0$, it follows from the definition of $r_1$ and (3.33) that

$$u(x'_1) \leq (1 + 8\varepsilon')^{\hat{y}_M} u(x_2) \leq (1 + 8\varepsilon')^{\hat{y}/M} \left(\frac{1+\varepsilon}{1-\varepsilon_0}\right)(\overline{G} \circledast u)(x_1).$$

So, since $\varepsilon < \varepsilon_0$, $\varepsilon' < \varepsilon_1$ and $M > 100$,

$$(3.35) \quad \begin{aligned} \log\left(\frac{(\overline{G} \circledast u)(x_1)}{u(x'_1)}\right) & \\ &\geq -\log\left(\frac{1+\varepsilon}{1-\varepsilon_0}\right) - \frac{\hat{y}}{M}\log(1+8\varepsilon') \\ &\geq -6(\varepsilon_1 \hat{y} + \varepsilon_0). \end{aligned}$$

On the other hand, by (3.3) and (3.4),

$$(3.36) \quad \begin{aligned} \log_b&\left(\frac{1 + \varepsilon_1 - u(x'_2)/u(x'_1)}{1 + \varepsilon_1 - (\overline{G} \circledast u)(x_2)/(\overline{G} \circledast u)(x_1)}\right) \\ &\geq \log_b\left(\frac{\varepsilon_1 - \varepsilon' + \delta e^{a\hat{y}/4}}{\varepsilon_1 - \varepsilon}\right) \geq \frac{a\hat{y}/4 + \log(\delta/(\varepsilon_1 - \varepsilon))}{\log b}. \end{aligned}$$

Together, the bounds in (3.35) and (3.36) imply (3.6). □

**4. Examples: branching random walks and cover time for the regular binary tree.** In this section we demonstrate Theorems 1.1 and 1.2 on the tightness of the maximal displacement of branching random walks (BRW) and the tightness of the cover time for the regular binary tree. The demonstration of each result consists of verifying the Assumptions 2.1, 2.2 and 2.4 in an appropriate setting, and then applying Theorem 2.5. The term $F_n$ in Theorem 2.5 corresponds to the distribution function of the maximal displacement of the BRW in the first setting. For the cover time problem, the relationship is less immediate and requires certain comparisons. In both settings the functions $Q_n$ and $\tilde{Q}$ that are employed in Assumptions 2.1 and 2.4 will be strictly concave, and will satisfy $\tilde{Q} = Q_n = Q$. As mentioned after Definition 2.3, properties (T1) and (T2) are therefore automatic.

4.1. *Branching random walk.* As in the introduction, we consider BRW whose underlying branching processes have offspring distribution $\{p_k\}_{k \geq 1}$ and whose random walk increments have distribution function $G$. As in Theorem 1.1, it is assumed that $\{p_k\}_{k \geq 1}$ satisfies $p_1 < 1$ and the moment condition (1.17), and that $G$ satisfies (G2) of Assumption 2.2. As before, we denote by $\mathcal{M}_n$ the maximal displacement of the BRW at time $n$, and by $F_n$ its distribution function. In the introduction we saw that $\{F_n\}_{n \geq 0}$ satisfies the recursion (1.3), with $Q(u) = 1 - \sum_{k=1}^{\infty} p_k(1-u)^k$. It is not difficult to check that the recursion (1.3) is a special case of the recursion given by (2.3), with $Q_n(\cdot) = Q(\cdot)$ and $G_n(y, \cdot) = G(\cdot)$ for all $n$ and $y$. This simplifies the checking required in Assumptions 2.1, 2.2 and 2.4.

PROOF OF THEOREM 1.1. Because of (G3), Assumption 2.2 need not be satisfied for the above choice of $G$. We therefore instead consider the translated BRW with increments having distribution function $G^{(L)}(x) = G(x - L)$, where $L > 0$ is chosen large enough to satisfy $G(-L) \leq \varepsilon_0 = (\log m)/100$, with $m = (1 + (m_1 - 1)/2) \wedge 2$. The maximal displacement at time $n$ of the corresponding BRW has distribution function $F_n^{(L)}(x)$ which satisfies the same recursion as $F_n$, with $G^{(L)}$ in place of $G$. One has $F_n^{(L)}(x) = F_n(x - nL)$ and, therefore, $(F_n^{(L)})^s = F_n^s$. So, tightness of $\{F_n^s\}_{n \geq 0}$ will follow from that of $(F_n^{(L)})^s$.

Theorem 1.1 will follow once we verify the assumptions of Theorem 2.5 for the BRW with increments distributed according to $G^{(L)}$. We first note that all of the conditions in Assumption 2.2 are satisfied. The condition (G1) for $G^{(L)}$ is immediate, and (G3) is satisfied because of our choice of $L$. The condition (G2) for $G$ implies (G2') for $G^{(L)}$, with $a' = a$ and $M' = 2L \vee M_0$, which in turn implies (G2) for $G^{(L)}$, with a new choice of $a$.

It is easy to see that (2.4) of Assumption 2.1 holds in a small enough neighborhood of 0, since $m < m_1$; we chose $m$ so that $m \in (1, 2]$. The bound (2.5) requires a little estimation. We write $Q(x) = m_1(x - g(x))$; one has $g(x) \in (0, x)$, and so for $0 < x_1 < x_2 \leq 2(\delta_0 \wedge x_1)$, with $\delta_0 > 0$ chosen small enough, it follows that

$$
\text{(4.1)} \quad \frac{Q(x_2)}{Q(x_1)} = \frac{x_2}{x_1} \left( \frac{1 - g(x_2)/x_2}{1 - g(x_1)/x_1} \right) \geq \frac{x_2}{x_1} \left( 1 - \frac{g(x_2)}{x_2} \right)
$$
$$
\geq \frac{x_2}{x_1} \frac{1}{1 + 2g(x_2)/x_2}.
$$

On the other hand, since $m_\theta < \infty$ for a given $\theta \in (1, 2]$, one has

$$
\text{(4.2)} \quad Q(x) \geq m_1 x - cx^\theta
$$

for appropriate $c > 0$ and for $x > 0$ close enough to 0 (see [25], page 212). It follows from (4.1) and (4.2) that, for small enough $\delta_0$,

$$
\frac{Q(x_2)}{Q(x_1)} \geq \frac{x_2}{x_1} \frac{1}{1 + 4c(Q(x_2))^{\theta - 1}}.
$$

Since $Q$ is concave and $Q(x_2) \leq 2Q(x_1)$, (2.5) follows from this, with $\theta^* = \theta - 1$ and $c^* = 8c$.

For Assumption 2.4, we only need to verify (2.10) and (2.11), with $\tilde{Q} = Q$ and $T_n u = -\overline{G}^{(L)} * Q(u)$. For $u \in \mathcal{D}$,

$$-(\overline{G}^{(L)} * Q(u))(x) \geq \overline{G}^{(L)}(-B)Q(u(x+B)) \geq Q(u(x+B)) - G^{(L)}(-B),$$

which implies (2.10), if $B = B(\eta_1)$ is chosen large enough so that $G^{(L)}(-B) < \eta_1$. On the other hand, for $B \geq 2L$,

$$-(\overline{G}^{(L)} * Q(u))(x) \leq Q(u(x-B)) - \int_{-\infty}^{x-B} \overline{G}^{(L)}(x-y)\, dQ(u(y))$$

$$\leq Q(u(x-B)) + e^{-aB/3},$$

where the last inequality follows from (G2'). Choosing $e^{-aB/3} < \eta_1$ implies (2.11).

We have demonstrated Assumptions 2.1, 2.2 and 2.4. It therefore follows from Theorem 2.5 that $\{\mathcal{M}_n^s\}_{n \geq 1}$ is tight, which implies Theorem 1.1. □

REMARK 4.1. As mentioned in the introduction, one can study the BRW satisfying the recursion

$$(4.3) \qquad \overline{F}_{n+1}^r = Q(G \circledast \overline{F}_n^r) = Q(-\overline{G} * \overline{F}_n^r),$$

rather than (1.3). In this setting the analog of Theorem 1.1 will continue to hold. To see this, note that since $\overline{F}_1 = -\overline{G} * \overline{F}_0$, one has $\overline{F}_n^r = Q(\overline{F}_n)$ for all $n$. Since $Q$ is continuous, with $Q(0) = 0$ and $Q(1) = 1$, the tightness of $\{(F_n^r)^s\}_{n \geq 0}$ follows directly from the tightness of $\{F_n^s\}_{n \geq 0}$.

REMARK 4.2. It is easy to see that if in (G2), $G_n(y, \cdot) = G(\cdot)$ for all $y$, with

$$(4.4) \qquad \overline{G}(x) = 0 \qquad \text{for } x \geq B_1$$

and some $B_1$, then (G2) holds, and so such $G$ provide a particular case of a BRW that is covered by Theorem 1.1. For such $G$, there are simple direct proofs of Theorem 1.1; see, for instance, [13, 27] and [11].

4.2. *Cover time for the regular binary tree.* In Section 1 we introduced the cover time $\mathcal{C}_n$ for the regular binary tree $\mathbf{T}_n$ of depth $n$, and in Theorem 1.2, we claimed that the sequence $\{\mathcal{E}_n^s\}_{n \geq 0}$ is tight, for $\mathcal{E}_n = \sqrt{\mathcal{C}_n/2^n}$. In this subsection we employ Theorem 2.5 to prove Theorem 1.2. As mentioned in Section 1, we will rely on work in Aldous [4], where it was shown that

$$\mathcal{C}_n/4(\log 2)n^2 2^n \to_{n \to \infty} 1 \qquad \text{in probability.}$$

It will be more convenient to instead demonstrate the tightness of $\{\tilde{\mathcal{E}}_n^s\}_{n \geq 0}$, where $\tilde{\mathcal{E}}_n = \sqrt{\tilde{\mathcal{C}}_n/2^n}$, and $\tilde{\mathcal{C}}_n$ is the cover time of the extended tree $\tilde{\mathbf{T}}_n$ that is formed from $\mathbf{T}_n$ by inserting an additional vertex $oo$ that is connected only to the root $o$ of $\mathbf{T}_n$. One can then deduce the tightness of $\{\mathcal{E}_n^s\}_{n\geq 0}$ from that of $\{\tilde{\mathcal{E}}_n^s\}_{n\geq 0}$.

For an appropriate sequence of random variables $\{\beta_n\}_{n\geq 0}$ related to $\{\tilde{\mathcal{E}}_n\}_{n\geq 0}$, Aldous showed that

$$\beta_{n+1} \stackrel{d}{=} K(\beta_n \vee \beta_n'), \tag{4.5}$$

where $\beta_0 = 0$, $\beta_n'$ is an independent copy of $\beta_n$, and, for any nonnegative random variable $Y$, $K(Y)$ denotes a nonnegative random variable that will be defined shortly which, when conditioned on $Y = y \geq 0$, will have a density $k(y, \cdot)$. If one sets $F_n(x) = P(\beta_n \leq x)$ and

$$\overline{G}(y, x-y) = \begin{cases} \int_x^\infty k(y,z)\,dz, & \text{for } y \geq 0, \\ \overline{G}(0,x), & \text{for } y < 0, \end{cases} \tag{4.6}$$

it follows from (4.5) and (4.6) that

$$\overline{F}_{n+1}(x) = (\overline{G} \circledast Q(\overline{F}_n))(x), \tag{4.7}$$

where $Q(u) = 2u - u^2$ is the concave mapping we have referred to repeatedly. [We note that in (4.7) the choice of $\overline{G}(y, x-y)$ for $y < 0$ is somewhat arbitrary, since $\beta_n \geq 0$)]. We will apply Theorem 2.5 to (4.7), and use this to deduce the desired tightness of $\{\tilde{\mathcal{E}}_n^s\}_{n\geq 0}$, and, hence, of $\{\mathcal{E}_n^s\}_{n\geq 0}$.

The sequence $\{\beta_n\}_{n\geq 0}$ that is employed in (4.5) is defined by $\beta_n = \sqrt{\Gamma(R_n)}$. Here, $R_n$ denotes the number of times the directed edge $(o, oo)$ is traversed by a nearest neighbor symmetric random walk by time $\tilde{\mathcal{C}}_n$, with $R_0 = 0$. For any nonnegative integer valued random variable $N$, $\Gamma(N)$ denotes a random variable which, when conditioned on $N = k$, has a Gamma$(k,1)$ distribution for $k \geq 1$, that is, it has a density $h(y) = y^{k-1}e^{-y}/(k-1)!$; we set $\Gamma(0) = 0$. The density $k(y, \cdot)$ employed above is

$$k(y,x) = \mathbf{1}_{x>0} 2x e^{-(x^2+y^2)} \sum_{j=0}^\infty \frac{(xy)^{2j}}{(j!)^2} \qquad \text{for } y \geq 0. \tag{4.8}$$

From (4.6) and (4.8), it follows that, for $y \geq 0$,

$$\overline{G}(y, x-y) = \begin{cases} e^{-(x^2+y^2)} \sum_{k=0}^\infty \frac{x^{2k}}{k!} \sum_{j=k}^\infty \frac{y^{2j}}{j!}, & x > 0, \\ 1, & x \leq 0. \end{cases} \tag{4.9}$$

For $k$ as in (4.8), Aldous showed that as $y \to \infty$, $\overline{G}(y, \cdot - y)$ converges to the normal distribution with mean 0 and variance $1/2$ (see [4], (9)). Since $\beta_n \to \infty$ as $n \to \infty$, this and (4.5) say heuristically that $\{\beta_n\}_{n\geq 0}$ should behave like the maximal displacement of a binary branching random walk with normal increments, which is a special case of the processes analyzed in Section 4.1.

Without going into details on Aldous' computations for the derivation of (4.5), we provide the following motivation. The random variable $\Gamma(R_n)$ corresponds to the cumulative time spent going from $o$ to $oo$ for the natural continuous time analog of the symmetric simple random walk, with transition rates equal to 1. The sequence $\{\Gamma(R_n)\}_{n\geq 0}$ will be easier to analyze than $\{R_n\}_{n\geq 0}$ itself. In particular, Aldous showed that

$$\Gamma(R_{n+1}) \stackrel{d}{=} \Gamma''(1 + \mathcal{P}(\Gamma(R_n) \vee \Gamma'(R_n'))), \tag{4.10}$$

where $R_n'$ is an independent copy of $R_n$, $\Gamma'(R_n')$ is a copy of $\Gamma(R_n')$ which is independent of $\Gamma(R_n)$ and, conditioned on $\Gamma(R_n) \vee \Gamma'(R_n') = y$, $\mathcal{P}(\cdot)$ is Poisson distributed with mean $y$ and $\Gamma''(1 + \mathcal{P}(\cdot))$ is a Gamma random variable of parameter $1 + \mathcal{P}(y)$. For $\beta_n = \sqrt{\Gamma(R_n)}$, it is not difficult to show from (4.10) that the conditional density $k(y,\cdot)$ of $K(Y)$ satisfies (4.8).

We will later derive the tightness of $\{\tilde{\mathcal{E}}_n^s\}_{n\geq 0}$ from the tightness of $\{\beta_n^s\}_{n\geq 0}$, and the tightness of $\{\mathcal{E}_n^s\}_{n\geq 0}$ from that of $\{\tilde{\mathcal{E}}_n^s\}_{n\geq 0}$. For the former, we note that the cover time $\tilde{\mathcal{C}}_n$ is sandwiched between the sum of $R_n$ and $R_n + 1$ random variables that correspond to the incremental return times to $oo$ of the random walk. These epochs will be i.i.d. and, when scaled by $2^n$, will have uniformly bounded variances. The fluctuations of $R_n$ will be at most of order $\sqrt{R_n}$ (since $\{\beta_n^s\}_{n\geq 0}$ is tight), and so one can show that the fluctuations of $\tilde{\mathcal{C}}_n/2^n$ will be at most of order $\sqrt{\tilde{\mathcal{C}}_n/2^n}$, which in turn will imply that the sequence $\{\tilde{\mathcal{E}}_n^s\}_{n\geq 0}$ is tight.

We now lay the groundwork for proving Theorem 1.2. We first demonstrate the following result.

THEOREM 4.3. *The sequence of random variables $\{\beta_n^s\}_{n\geq 0}$ for the regular binary tree is tight.*

In order to apply Theorem 2.5 to the sequence $\{\beta_n^s\}_{n\geq 0}$, we need to verify that Assumptions 2.1, 2.2 and 2.4 hold. The choice of $Q(u) = 2u - u^2$ here is a special case of that considered in Section 4.1, so Assumption 2.1 holds for this $Q$. As in Section 4.1, $\overline{G}(y, x-y)$ need not satisfy (G3) of Assumption 2.2. We handle this in a way similar to what was done there, by translating $\overline{G}$. Set $\varepsilon_0 = (\log m)/100$, with $m$ as in Assumption 2.1, and let $\mathcal{N}(0, 1/2)$ denote a mean zero Gaussian random variable of variance $1/2$. We fix a constant $L > 1$ such that

$$P(\mathcal{N}(0, 1/2) > -L) \geq 1 - \varepsilon_0, \tag{4.11}$$

and set $\beta_n^{(L)} = \beta_n + nL$. Defining $\overline{F}_n^{(L)} = P(\beta_n^{(L)} > x)$, it follows from (4.7) that

$$\overline{F}_{n+1}^{(L)}(x) = (\overline{G}_n^{(L)} \circledast Q(\overline{F}_n^{(L)}))(x), \tag{4.12}$$

where

$$\overline{G}_n^{(L)}(y, x) = \overline{G}(y - nL, x - L).$$

Note that the tightness of $\{(F_n^{(L)})^s\}_{n \geq 0}$ is equivalent to the tightness of $\{F_n^s\}_{n \geq 0}$. We proceed to verify that Assumptions 2.1, 2.2 and 2.4 hold for the recursions (4.12).

As before, Assumption 2.1 holds with the function $Q(u) = 2u - u^2$. In order to verify Assumption 2.2, we first need to verify (G1) for $\overline{G}_n^{(L)}$; this follows immediately from the analogous statement on $\overline{G}$:

LEMMA 4.4. *The function $\overline{G}(y, x - y)$ is increasing in $y$, whereas $\overline{G}(y, x)$ is decreasing in $y$.*

PROOF OF LEMMA 4.4. The claim is obvious when $y < 0$. For $y \geq 0$, making the transformation $y^2 \mapsto s$, one has

$$\frac{\partial \overline{G}(y, x - y)}{\partial y} = 2y \frac{\partial \overline{G}(\sqrt{s}, x - \sqrt{s})}{\partial s} \bigg|_{s = y^2}.$$

On the other hand, from (4.9),

$$\begin{aligned}
\frac{\partial \overline{G}(\sqrt{s}, x - \sqrt{s})}{\partial s} \\
&= e^{-s} e^{-x^2} \left( -\sum_{j=0}^{\infty} \frac{s^j}{j!} \sum_{k=0}^{j} \frac{x^{2k}}{k!} + \sum_{j=0}^{\infty} \frac{j s^{j-1}}{j!} \sum_{k=0}^{j} \frac{x^{2k}}{k!} \right) \\
&= e^{-s} e^{-x^2} \left( -\sum_{j=0}^{\infty} \frac{s^j}{j!} \sum_{k=0}^{j} \frac{x^{2k}}{k!} + \sum_{j=0}^{\infty} \frac{s^j}{j!} \sum_{k=0}^{j+1} \frac{x^{2k}}{k!} \right) \\
&= e^{-s} e^{-x^2} \left( \sum_{j=0}^{\infty} \frac{s^j}{j!} \frac{x^{2(j+1)}}{(j+1)!} \right) > 0.
\end{aligned} \tag{4.13}$$

This completes the proof of the required monotonicity of $\overline{G}(y, x - y)$ in $y$.

The proof of the required monotonicity of $\overline{G}(y, x)$ in $y$ is more subtle. First, substituting in (4.9), one has, for $x > -y$, that

$$\overline{G}(y, x) = e^{-(y^2 + (x+y)^2)} \sum_{j=0}^{\infty} \frac{y^{2j}}{j!} \sum_{k=0}^{j} \frac{(x+y)^{2k}}{k!}.$$

Thus, after some algebra,

$$\frac{1}{2}e^{y^2+(x+y)^2}\frac{\partial \overline{G}(y,x)}{\partial y}$$

$$= \left[-y\sum_{j=0}^{\infty}\frac{y^{2j}}{j!}\sum_{k=0}^{j}\frac{(x+y)^{2k}}{k!} + \sum_{j=0}^{\infty}\frac{jy^{2j-1}}{j!}\sum_{k=0}^{j}\frac{(x+y)^{2k}}{k!}\right.$$

$$\left. -(x+y)\sum_{j=0}^{\infty}\frac{y^{2j}}{j!}\sum_{k=0}^{j}\frac{(x+y)^{2k}}{k!} + \sum_{j=0}^{\infty}\frac{y^{2j}}{j!}\sum_{k=0}^{j}\frac{k(x+y)^{2k-1}}{k!}\right]$$

(4.14)

$$= \left[-y\left(\sum_{j=0}^{\infty}\frac{y^{2j}}{j!}\sum_{k=0}^{j}\frac{(x+y)^{2k}}{k!} - \sum_{j=0}^{\infty}\frac{y^{2j}}{j!}\sum_{k=0}^{j+1}\frac{(x+y)^{2k}}{k!}\right)\right.$$

$$\left. -(x+y)\left(\sum_{j=0}^{\infty}\frac{y^{2j}}{j!}\sum_{k=0}^{j}\frac{(x+y)^{2k}}{k!} - \sum_{j=1}^{\infty}\frac{y^{2j}}{j!}\sum_{k=0}^{j-1}\frac{(x+y)^{2k}}{k!}\right)\right]$$

$$= (x+y)\sum_{j=0}^{\infty}\frac{(y(x+y))^{2j}}{(j!)^2}\left(\frac{y(x+y)}{j+1}-1\right) = (x+y)g(y(x+y))$$

where

$$g(z) = \sum_{j=0}^{\infty}\frac{z^{2j}}{(j!)^2}\left(\frac{z}{j+1}-1\right).$$

The required monotonicity of $\overline{G}(y,x)$ therefore follows once we prove that $g(z/2)$ is nonpositive for $z \geq 0$. But, for $I_\nu(x)$ denoting the modified Bessel function of order $\nu$ (see, e.g., [1], (9.6.10) for the definition), one has

$$g(z/2) = I_1(z) - I_0(z) = -\frac{1}{\pi}\int_0^\pi e^{z\cos\theta}(1-\cos\theta)\,d\theta \leq 0,$$

where the equality is a consequence of [1], (9.6.19). This completes the proof of Lemma 4.4. $\square$

We still need to check (G2) and (G3). The latter is immediate from our choice of $L$ and the following fact from [4], equation (9):

(4.15) $$\lim_{y\to\infty}\overline{G}(y,x) = P(\mathcal{N}(0,1/2) > x).$$

In particular, by Lemma 4.4, $\overline{G}(y,x)$ decreases in $y$ and, therefore,

(4.16) $$\overline{G}_n^{(L)}(y,0) = \overline{G}(y-nL,-L) \geq \lim_{y\to\infty}\overline{G}(y,-L) \geq 1-\varepsilon_0,$$

where (4.15) and (4.11) were used in the last step.

It thus only remains to verify that, (G2) also holds for $\overline{G}_n^{(L)}$. Because (G2′) implies (G2), it suffices to verify that, for $M = 4L$, there exists an $a > 0$ such that, for all $x \geq 0$,

$$(4.17) \qquad \overline{G}(y - M, x + M) \leq e^{-aM}\overline{G}(y, x).$$

It is clearly enough to consider only the case $y \geq M$, since $\overline{G}(y, x)$ is decreasing in $y$ and $\overline{G}(y, x) = \overline{G}(0, x)$ for $y \leq 0$. The inequality (4.17) therefore follows from the following lemma.

LEMMA 4.5. *For $M$ and $L$ as above, there exists $C = C(L) > 0$ such that, for all $y > M/2$ and $x \geq y$,*

$$(4.18) \qquad \frac{\partial \overline{G}(y, x - y)/\partial y}{\overline{G}(y, x - y)} \geq C.$$

PROOF OF LEMMA 4.5. Throughout this proof, $C_i = C_i(L)$ denote strictly positive constants. From (4.13), it follows that

$$\frac{\partial \overline{G}(y, x - y)/\partial y}{\overline{G}(y, x - y)} = \frac{xye^{-(x^2+y^2)}I_1(2xy)}{\overline{G}(y, x - y)},$$

where we recall that $I_1(z)$ is the modified Bessel function of order 1. Since $I_1(z)$ is asymptotic to $e^z/\sqrt{2\pi z}$ for $z$ large (see [1], (9.7.1)), and is positive and continuous for $z > 0$, there exists a constant $C_2 > 0$ such that $\sqrt{xy}e^{-2xy}I_1(2xy) \geq C_2$ for $xy > 1$. The lemma thus follows once we show that there exists $C_3 > 0$ such that, for $y > 2$ and $x > y$,

$$(4.19) \qquad \overline{G}(y, x - y) \leq C_3\sqrt{xy}e^{-(x-y)^2}.$$

For large enough $C_3$, (4.19) obviously holds when $|x - y| \leq 1$. (Recall that $x, y > 2$.) Therefore, it suffices to prove (4.19) for $x > y + 1$. To show this, rewrite (4.9) for $x, y > 2$ as

$$\overline{G}(y, x - y) = e^{-(x^2+y^2)}\left[\sum_{k=0}^{\lfloor y^2 \rfloor} + \sum_{k=\lfloor y^2 \rfloor + 1}^{\infty}\right]\frac{x^{2k}}{k!}\sum_{j=k}^{\infty}\frac{y^{2j}}{j!},$$

which we rewrite as $G_1(y, x - y) + G_2(y, x - y)$. First, note that by replacing the summation over $j \in \{k, k+1, \ldots\}$ by a summation over $j \in \mathbb{N}$, one gets, using Stirling's formula,

$$G_1(y, x - y) \leq e^{-x^2}\sum_{k=0}^{\lfloor y^2 \rfloor}\frac{x^{2k}}{k!} \leq e^{-x^2}\sum_{k=0}^{\lfloor y^2 \rfloor}\left(\frac{ex^2}{k}\right)^k$$

$$\leq C_4 e^{-x^2}\left(\frac{ex^2}{y^2}\right)^{y^2},$$

for some $C_4 > 0$. But, for $x > y > 0$, since $\log u \leq u - 1$ for $u > 0$,
$$\log(ex^2/y^2) \leq \frac{2x}{y} - 1$$
and, therefore,
$$(4.20) \qquad G_1(y, x - y) \leq C_4 e^{-(x-y)^2}.$$

To handle $G_2(y, x - y)$, we note that in the summation there, $j > y^2$ and, hence, again applying Stirling's formula,
$$\sum_{j=k}^{\infty} \frac{y^{2j}}{j!} \leq C_5 \sqrt{k} \left(\frac{ey^2}{k}\right)^k.$$

Therefore, another application of Stirling's formula yields
$$G_2(y, x - y) \leq C_6 e^{-(x^2 + y^2)} \sum_{k=\lfloor y^2 \rfloor}^{\infty} \left(\frac{exy}{k}\right)^{2k}$$
$$\leq C_7 \sqrt{xy} e^{-(x-y)^2}.$$

Together with (4.20), this completes the proof of (4.19) and, hence, of Lemma 4.5. $\square$

We still need to verify Assumption 2.4 for the kernels $\overline{G}_n^{(L)}$ and the transformation $(T_n u)(x) = (\overline{G}_n^{(L)} \circledast u)(x)$. For $u \in \mathcal{D}$,
$$(T_n u)(x) \geq \overline{G}_n^{(L)}(x + B, -B) Q(u(x + B)) \geq Q(u(x + B)) - G_n^{(L)}(x + B, -B).$$
But by (G1),
$$G_n^{(L)}(x + B, -B) \leq \lim_{y \to \infty} G(y, -B) = P(\mathcal{N}(0, 1/2) \leq -B),$$
which implies (2.10) if one chooses $B = B(\eta_1)$ such that $P(\mathcal{N}(0, 1/2) \leq -B) < \eta_1$. Similarly, taking $B = B(\eta_1)$ to be a multiple of $M$ and using (G2'),
$$T_n u(x) \leq Q(u(x - B)) - \int_{-\infty}^{x-B} \overline{G}^{(L)}(y, x - y) \, dQ(u(y))$$
$$\leq Q(u(x - B)) + e^{-aB/3},$$
which implies (2.11) for $e^{-aB/3} < \eta_1$.

We have verified Assumptions 2.1, 2.2 and 2.4 for the recursions in (4.12). Theorem 2.5 therefore holds for $\{F_n^{(L)}\}_{n \geq 0}$, and therefore for $\{F_n\}_{n \geq 0}$. We thus obtain Theorem 4.3.

As outlined in the beginning of this subsection, one can use Theorem 4.3 to prove Theorem 1.2.

PROOF OF THEOREM 1.2. First note that by coupling the random walks on $\mathbf{T}_n$ and $\tilde{\mathbf{T}}_n$ in the natural way, one has $\mathcal{C}_n \leq \tilde{\mathcal{C}}_n \leq \mathcal{C}_n + 2R_n$. From [4], (8) it follows that $ER_n < 2n^2 \log 2$ and also that, for $\delta > 0$ small enough, $P(\mathcal{C}_n \leq \delta^2 2^n) \to_{n \to \infty} 0$. So, for such $\delta$,

$$\lim_{n \to \infty} P\left(\sqrt{\frac{\tilde{\mathcal{C}}_n}{2^n}} - \sqrt{\frac{\mathcal{C}_n}{2^n}} > \delta\right) \leq \lim_{n \to \infty} P(2R_n \geq \delta^2 2^n)$$

$$\leq \lim_{n \to \infty} \frac{4n^2 \log 2}{\delta^2 2^n} = 0.$$

Hence, the tightness of $\{\mathcal{E}_n^s\}_{n \geq 0}$ follows if one proves the analogous result for $\tilde{\mathcal{C}}_n$ instead of $\mathcal{C}_n$.

Letting $\tau_i^{(n)}$ denote the length of the $i$th epoch between returns to $oo$ of the random walk in the extended tree of depth $n$, one has

$$(4.21) \qquad \sum_{i=1}^{R_n} \tau_i^{(n)} \leq \tilde{\mathcal{C}}_n \leq \sum_{i=1}^{R_n+1} \tau_i^{(n)}.$$

In the proof we will employ the limits

$$(4.22) \qquad \lim_{n \to \infty} \frac{E\tau_1^{(n)}}{2^{n+1}} = 2, \qquad \lim_{n \to \infty} \frac{\operatorname{Var}(\tau_1^{(n)})}{2^{2(n+1)}} = 12.$$

To see (4.22), note that $\tau_1^{(n)}$ is identical in law to the return time to 0 of a simple random walk $W_j$ on $\{0, 1, \ldots, n+1\}$, having probability of jump to the right equal to $2/3$ and to the left equal to $1/3$. This return time is asymptotically equivalent to the product of a Bernoulli($1/2$) random variable with the sum of a geometric number of i.i.d. random variables, each corresponding to the time it takes for the walk $W_j$ started at $n$ to return to $n$ after hitting either 0 or $n+1$; the parameter of the geometric random variable is the probability that such a random walk hits 0 before hitting $n+1$, and all random variables involved are independent. Standard computations, using, for example, [16], page 314 (2.4) and page 317 (3.4) lead to (4.22).

Since $\beta_n = \sqrt{\Gamma(R_n)}$, it follows from Theorem 4.3 that, for an appropriate deterministic sequence $\{B_n\}$, $B_n \geq 1$, there exist constants $\varepsilon_J \to_{J \to \infty} 0$ such that

$$(4.23) \qquad P(|\sqrt{\Gamma(R_n)} - B_n| > J) \leq \varepsilon_J.$$

One can check that

$$(4.24) \qquad P(|\Gamma(R_n) - B_n^2| > 3J^2 B_n) \leq \varepsilon_J.$$

The random variable $\Gamma(R_n)$ is the sum of $R_n$ i.i.d. exponentials that, conditioned on $R_n$, is asymptotically Gaussian with variance $R_n$. There therefore exists a deterministic sequence $\{A_n\}$ such that, setting

$$C_K(n) = P(|R_n - A_n| \geq K\sqrt{A_n}),$$

one has

(4.25) $$\lim_{K\to\infty}\limsup_{n\to\infty} C_K(n)=0.$$

Setting $\bar\tau_i^{(n)}=\tau_i^{(n)}/E\tau_1^{(n)}$, one gets, using (4.22) in the first inequality and the inequality $\sqrt{1+x}-1\le x/2$ in the second, that, for any $K_1>0$,

$$\lim_{n\to\infty} P\left(\left|\sqrt{\frac{\tilde{\mathcal{C}}_n}{2^n}}-\sqrt{\frac{A_n E\tau_1^{(n)}}{2^n}}\right|\ge K_1\right)\le \lim_{n\to\infty} P\left(\left|\sqrt{\frac{\tilde{\mathcal{C}}_n}{E\tau_1^{(n)}}}-\sqrt{A_n}\right|\ge K_1\right)$$

$$\le \lim_{n\to\infty} P\left(\left|\frac{\tilde{\mathcal{C}}_n}{E\tau_1^{(n)}}-A_n\right|\ge K_1\sqrt{A_n}\right).$$

By (4.21), on the event $\{|R_n-A_n|\le K_1\sqrt{A_n}/4\}$,

$$\sum_{i=1}^{A_n-K_1\sqrt{A_n}/4}\tau_i^{(n)}\le \tilde{\mathcal{C}}_n\le \sum_{i=1}^{A_n+K_1\sqrt{A_n}/4+1}\tau_i^{(n)}.$$

Since $A_n\to\infty$, using this decomposition in the first inequality, and the fact that the $\tau_i^{(n)}$ are i.i.d, together with (4.22) and Markov's inequality in the second, one obtains

$$\lim_{n\to\infty} P\left(\left|\sqrt{\frac{\tilde{\mathcal{C}}_n}{2^n}}-\sqrt{\frac{A_n E\tau_1^{(n)}}{2^n}}\right|\ge K_1\right)$$

(4.26)
$$\le \lim_{n\to\infty} P\left(|R_n-A_n|\ge \frac{K_1\sqrt{A_n}}{4}\right)$$
$$+\lim_{n\to\infty} P\left(\sum_{i=1}^{A_n+K_1\sqrt{A_n}/4+1}(\bar\tau_i^{(n)}-1)\ge \frac{3K_1\sqrt{A_n}}{4}\right)$$
$$+\lim_{n\to\infty} P\left(\sum_{i=1}^{A_n-K_1\sqrt{A_n}/4}(\bar\tau_i^{(n)}-1)\le -\frac{3K_1\sqrt{A_n}}{4}\right)$$
$$\le \lim_{n\to\infty} C_{K_1/4}(n)+\frac{8}{K_1^2}.$$

Together with (4.25), this demonstrates the tightness of $\{\mathcal{E}_n^s\}_{n\ge 0}$.

To demonstrate (1.20), we argue by contradiction. Assuming that (1.20) does not hold, one can use the analog of the argument from 4.23 to 4.24 [with $\mathcal{C}_n/2^n$ replacing $\Gamma(R_n)$] and the fact that $\mathcal{C}_n/n^2 2^n \to 4\log 2$ to deduce that, for any $\varepsilon_2>0$, there exist constants $b_n<c_n$, with $c_n-b_n=\varepsilon_2 n 2^n/100$, such that, for large enough $n$,

(4.27) $$P^o(\mathcal{C}_n\in(b_n,c_n))\ge \tfrac{7}{8}.$$

Here and in the remainder of the proof, we use the notation $P^v$ to denote the law of the symmetric simple random walk $\{X_j\}_{j\geq 0}$, started at a vertex $v$, on the (nonextended) tree $\mathbf{T}_n$. In particular,

$$P^o(\mathcal{C}_n \leq b_n) \leq 1/8. \tag{4.28}$$

Define $a_n = b_n - 1000 \cdot 2^n$, noting that $a_n > 0$ because $b_n/n^2 2^n \to 4\log 2$. We will show that, with high probability, $X_j = o$ at some $j \in (a_n, b_n)$. To see this, let $\rho_i^{(n)}$, $i \geq 1$, with $\rho_0^{(n)} = 0$, denote the successive return times to $o$ of the random walk $\{X_j\}_{j\geq 0}$ on $\mathbf{T}_n$, and set $\hat{\tau}_i^{(n)} = \rho_i^{(n)} - \rho_{i-1}^{(n)}$. For any $v \in \mathbf{L}_n$, one has by an application of the strong Markov property that

$$P^o(\{\rho_i^{(n)}\}_{i\geq 1} \cap (a_n, b_n) = \varnothing) = E(P^{X_{a_n-1}}(X_j \neq o, j = 1, \ldots, b_n - a_n))$$

$$\leq P^v(\hat{\tau}_1^{(n)} \geq b_n - a_n) \leq \frac{E^v(\hat{\tau}_1^{(n)})}{1000 \cdot 2^n}.$$

(Recall that for $n > 1$, $\mathbf{L}_n$ denotes the set of vertices at distance $n$ from $o$.) On the other hand,

$$E^o(\hat{\tau}_1^{(n)}) \geq P^o(T_{\mathbf{L}_n}^{(n)} < \hat{\tau}_1^{(n)}) E^v(\hat{\tau}_1^{(n)}) \geq \tfrac{1}{2} E^v(\hat{\tau}_1^{(n)}),$$

where $T_{\mathbf{L}_n}^{(n)} = \min\{j \geq 1 : X_j \in \mathbf{L}_n\}$. The first limit in (4.22) holds for the nonextended tree $\mathbf{T}_n$ as well, with $n$ replacing $n+1$. Together with the previous two displays, this implies that, for large $n$,

$$P^o(\{\rho_i^{(n)}\}_{i\geq 1} \cap (a_n, b_n) = \varnothing) \leq \frac{2E^o(\hat{\tau}_1^{(n)})}{1000 \cdot 2^n} \leq \frac{1}{200}. \tag{4.29}$$

Let $\rho_n^* = \inf\{\rho_i^{(n)} : \rho_i^{(n)} > a_n\}$ and $\mathcal{A}_n = \{\rho_n^* \in (a_n, b_n), \mathcal{C}_n > \rho_n^*\}$. It follows from (4.28) and (4.29) that

$$P^o(\mathcal{A}_n) \geq \tfrac{7}{8} - \tfrac{1}{200} \geq \tfrac{1}{2}. \tag{4.30}$$

On $\mathcal{A}_n$ there is at least one (random) vertex in $\mathbf{L}_n$ that has not been covered by the random walk by time $\rho_n^*$. Let $T_v^{(n)}$ denote the hitting time of $v \in \mathbf{T}_n$ after time 1. We will now show that one may chose $\varepsilon_2 > 0$ such that

$$P^o(T_v^{(n)} < c_n - a_n) \leq P^o(T_v^{(n)} < 2\varepsilon_2 n 2^n/100) < \tfrac{1}{2}. \tag{4.31}$$

The first inequality is immediate from $c_n - a_n < 2\varepsilon_2 n 2^n/100$. For the second inequality, first note that $P^o(T_v^{(n)} < \hat{\tau}_1^{(n)}) = 1/2n$. (One can see this by considering the symmetric simple random walk on the ray connecting $0$ and $v$, disregarding excursions off the ray.) By considering the set $\{T_v^{(n)} < \rho_{\lfloor n/4 \rfloor}^{(n)}\}$ and its complement separately, it follows from this that

$$P^o(T_v^{(n)} < 2\varepsilon_2 n 2^n/100) \leq \frac{n}{4} \cdot \frac{1}{2n} + P^o\left(\sum_{i=1}^{\lfloor n/4 \rfloor} \hat{\tau}_i^{(n)} < 2\varepsilon_2 n 2^n/100\right). \tag{4.32}$$

Note that by [16], page 314, (2.4), $P^v(\hat{\tau}_1^{(n)} < T_{\mathbf{L}_n}^{(n)}) = 1/(2^n - 1)$ and, hence, by the Markov property, for all $x > 0$,

$$P^o(\hat{\tau}_1^{(n)} > \lfloor x 2^n \rfloor) \geq \frac{1}{2} P^v(\hat{\tau}_1^{(n)} > \lfloor x 2^n \rfloor) \geq \frac{1}{2}\left(1 - \frac{1}{(2^n - 1)}\right)^{\lfloor x 2^n \rfloor}.$$

This bound implies that the random variables $\hat{\tau}_1^{(n)}/2^n$ possess exponential tails and have expectations bounded away from 0, uniformly in $n$. The second inequality in (4.31) follows from this and (4.32).

It follows from the strong Markov property and a little work that

$$P^o(\mathcal{A}_n, \mathcal{C}_n \geq c_n) \geq P^o(\mathcal{A}_n) \min_{\{v \in \mathbf{L}_n\}} P^o(T_v^{(n)} \geq c_n - a_n).$$

This is

$$\geq \frac{1}{2}\left(1 - \max_{\{v \in \mathbf{L}_n\}} P^o(T_v^{(n)} < c_n - a_n)\right) \geq \frac{1}{4},$$

where we have used (4.30) in the first inequality and (4.31) in the second. This bound contradicts (4.27). We have thus shown (1.20), which completes the proof of the theorem. □

REMARK 4.6. We have shown Theorem 1.2 for regular binary trees. The statement and its proof extend to regular $k$-ary trees. For general Galton–Watson trees with offspring distribution satisfying $p_0 = 0$, $p_1 < 1$ and (1.17), one can show the analog of Theorem 4.3 in the annealed setting, without in essence changing the proof. [One replaces the function $Q(u) = 2u - u^2$ by $Q(u) = 1 - \sum_{k=1}^{\infty} p_k(1-u)^k$.] The remainder of the proof of Theorem 1.2, however, relies on the regularity of the binary tree for various estimates, and so an extension to more general Galton–Watson trees would require additional effort.

**5. An open problem.** Consider the lattice torus $\mathbb{Z}_n^2 = \mathbb{Z}^2/n\mathbb{Z}^2$, and let $\mathbf{C}_n$ denote the number of steps required for a simple random walk to cover $\mathbb{Z}_n^2$. Confirming a conjecture in [3], it was proved in [14] that $\pi \mathbf{C}_n/4n^2(\log n)^2 \to 1$ in probability. The intuition (although not the details) behind the proof in [14] draws heavily from the covering of the regular binary tree by a simple random walk. One thus expects that a result similar to Theorem 1.2 should hold for $\mathbf{C}_n$. We therefore put forward the following conjecture.

CONJECTURE 5.1. *The sequence of random variables*

$$\mathbf{E}_n = \sqrt{\frac{\mathbf{C}_n}{n^2}} - \mathrm{Med}\left(\sqrt{\frac{\mathbf{C}_n}{n^2}}\right)$$

*is tight and nondegenerate.*

**Acknowledgment.** We thank the referees and Ming Fang for their careful reading of the paper and their useful suggestions.

SCHOOL OF MATHEMATICS
UNIVERSITY OF MINNESOTA
206 CHURCH ST. SE
MINNEAPOLIS, MINNESOTA 55455
USA
E-MAIL: bramson@math.umn.edu

SCHOOL OF MATHEMATICS
UNIVERSITY OF MINNESOTA
206 CHURCH ST. SE
MINNEAPOLIS, MINNESOTA 55455
USA
AND
FACULTY OF MATHEMATICS
WEIZMANN INSTITUTE OF SCIENCE
REHOVOT 76100
ISRAEL
E-MAIL: zeitouni@math.umn.edu